\documentclass[10pt,twoside]{article}
\usepackage{nevin-pf}
\usepackage{url}

\newcommand{\FF}{{}_2F_1}
\newcommand{\JF}[1]{}

\DeclareMathOperator{\Arg}{Arg}

\def\frak{\mathfrak}
\def\bb{\mathbb}
\def\cal{\mathcal}
\def\E{{\bb{E}}
\def\Pr{\bb P}}
\def\ds{\displaystyle}
\def\C{\bb C}
\def\Z{\bb Z}
\begin{document}

\title{\bf 
Singularity Analysis, \\ Hadamard Products, \\ and Tree Recurrences}
\author{James Allen Fill$^{1}$ \and Philippe Flajolet${}^2$ \and Nevin 
Kapur${}^1$\\  
\small ${}^1$  The Johns Hopkins University,
Department of Mathematical Sciences,\\[-2truept]
\small 34th and Charles Streets,
Baltimore, MD 21218-2682 (USA)\\[-2truept]
\small \texttt{\{jimfill,nevin\}@jhu.edu},
\texttt{http://www.mts.jhu.edu/\~{}\{fill,kapur\}/}\\
\small ${}^2$ {\sc Algorithms} Project, INRIA, F-78153 Rocquencourt
(France) \\[-2truept]
\small \texttt{Philippe.Flajolet@inria.fr},
\texttt{http://algo.inria.fr/flajolet/}%
}

\date{June~12, 2003}

\maketitle

\begin{abstract}
  We present a toolbox for extracting asymptotic information on the
  coefficients of combinatorial generating functions. This toolbox
  notably includes a treatment of the effect of Hadamard products on
  singularities in the context of the complex Tauberian technique
  known as singularity analysis. As a consequence, it becomes possible
  to unify the analysis of a number of divide-and-conquer algorithms,
  or equivalently random tree models, including several classical
  methods for sorting, searching, and dynamically managing equivalence
  relations.
\end{abstract}

\bigskip



This study was
motivated by a desire to unify the analysis of a number of algorithms
and data structures of computer science. By analysis
we mean 
here (precise) average-case analysis of cost functions as introduced
by Knuth and illustrated in the collection~\cite{Knuth00} as well as
in his monumental series, \emph{The Art of Computer Programming} (see
especially~\cite{Knuth97,Knuth98a}).  In the first part of this paper
(Section~\ref{DandC-sec} and~\ref{sec:some-special-tree}), we consider
a major paradigm of algorithmic design, the ``divide-and-conquer''
principle, which is closely related to families of random trees and
associated ``tree recurrences''.  The basic framework is described in
Section~\ref{DandC-sec}, while lead examples are introduced
in Section~\ref{sec:some-special-tree} below. Our treatment rests on
combinatorial generating functions.


The central part of this paper (Sections~\ref{sec:sing-expans-diff}
and~\ref{sec:hadam-prod-transf}) is devoted to the process of
extracting coefficients, at least asymptotically, from generating
functions.  Singularities have long been recognized to contain highly
useful information in this regard, and we start by recalling in
Section~\ref{sec:sing-expans-diff} the basic principles of the complex
Tauberian approach known as ``\emph{singularity analysis}''.
Applications to algorithms and trees require, in particular,
 techniques for coping with generating functions that may be
constructed by a tower of several transformations. Here, we develop
the theory of composition of singularities under Hadamard products in
Section~\ref{sec:hadam-prod-transf}.  (The reader only interested in
complex-analytic aspects can jump directly to
Sections~\ref{sec:sing-expans-diff} and~\ref{sec:hadam-prod-transf}.)

The final part (Sections~\ref{sec:some-applications}
and~\ref{sec:appl-high-moments}) returns to the original problem of
analysing divide-and-conquer algorithms, taking full advantage of the
analytic results of previous sections. Tree recurrences and first
moments form the subject of Section~\ref{sec:some-applications}, where
full  asymptotic expansions are derived for
expectations of costs.  Section~\ref{sec:appl-high-moments} describes
possible extensions of the basic framework to the determination of variances
and higher moments as well as to some other random tree models.

\section{Introduction}\label{DandC-sec}

``Divide-and-Conquer'' is a major principle of algorithmic design in
computer science. An instance ($I$) of a problem to be solved is first
split into smaller subproblems ($I',I''$) that are solved recursively
by the same process; the partial solutions are then woven back to
yield a solution to the original problem. The abstract scheme is then
of the form:
\begin{equation}\label{eq1}
\hbox{\begin{tabular}{lll}
solve$(I)$ &:=& $(I',I'')$ := split$(I)$;\\
        && $J'$ := solve$(I')$; $J''$ := solve$(I'')$;\\
        && return weave$(J',J'')$.
\end{tabular}}
\end{equation}
(Problems of size smaller than a certain threshold are treated
directly without any recursive call.)  Algorithms resorting to the
scheme~(\ref{eq1}) include classical sorting methods (mergesort,
quicksort, radix-exchange sort), data structures based on trees
(binary search trees, digital trees known as ``tries'', quadtrees for
multidimensional search, union--find trees) as well as various methods
used in computational geometry, distributed computation, and
communication theory. We refer the reader to classical books on data
structures, algorithms, and analysis of algorithms for details, for
instance,~\cite{CoLeRi90,GrKn82,Hofri95,Knuth98a,Mahmoud92,Mahmoud00,Sedgewick88,SeFl96,Szpankowski01,ViFl90}.

In general, a class of probabilistic models $\frak{M}_n$ indexed by
the size~$n$ of the problem instance is assumed to reflect the nature
of data fed to the algorithm.  A cost function---typically, the number
of certain operations performed by the algorithm---then becomes a
random variable~$X_n$ whose form is induced by~$\frak{M}_n$ and the
particular divide-and-conquer algorithm considered.  The problem is
then to obtain characteristics of~$X_n$, for instance its mean, higher
moments, or even distributional information.  The asymptotic
limit~$n\to\infty$ is usually considered, since an important
phenomenon of ``asymptotic simplification'' is to be expected in a
large number of situations.

Under natural conditions,
a recurrence that closely mimics the recursive structure
of~(\ref{eq1}) relates the random variables~$X_n$: 
\begin{equation}\label{eq2}
X_n=t_n+X_{K_n}+\widetilde{X}_{n-a-K_n}.
\end{equation} 
The interpretation is as follows: $t_n$ is a quantity\footnote{%
  Some analyses require a randomly varying toll. For mean value
  analysis, the distinction between deterministic and stochastic tolls
  is, however, immaterial.  }, called the ``toll'', that
  represents the 
cost incurred by splitting the initial instance and weaving back the
final solution; $K_n$ is the (random) size of the first subproblem, in
which case, the second subproblem has a size that is the complement of
$K_n$ to $n-a$, for some small constant $a$ (usually, $a=0$ or~$a=1$),
which is specific to the algorithm considered.  The random variables
of type $X$ and $K$ are assumed to be independent, as are
  the two $X$-sequences~$X$ and~$\widetilde{X}$ on the
  right  
  in~(\ref{eq2}), and a subproblem of
size~$k$ is assumed to satisfy model~$\mathfrak{M}_k$---this property
is sometimes called ``randomness preservation'' and is satisfied by
many cases of algorithmic interest.  A direct asymptotic treatment of
the recursive relation~(\ref{eq2}) binding random variables is
sometimes feasible; see the (metric) ``contraction method'' surveyed
by R\"osler and R\"uschendorf~\cite{RoRu01} and applied by
Neininger~\cite{Neininger02} to a subset of the problems discussed
here.

Turning to average-case analysis, the expected cost $f_n := \E(X_n)$
satisfies a \emph{recurrence} that is directly implied by~(\ref{eq2}):
\begin{equation}\label{eq3}
f_n=t_n+\sum_{k} p_{n,k} (f_k+f_{n-a-k}).
\end{equation}
with
the \emph{splitting probabilities} $p_{n,k}:=\Pr(K_n=k)$ being
determined by the model~$\frak{M}_n$ used.  Trees are naturally
associated with recursive procedures, and, accordingly, the
recurrence~(\ref{eq3}) can be viewed as associated with a random tree
model of the following form: the root has size~$a$, the left subtree
has size~$k$ with probability $p_{n,k}$, and the right subtree has the
remaining quantity~$n-a-k$ as size. Then~(\ref{eq3})is
interpreted as giving the expectation of a cost function over the tree
structure that is induced by the family of tolls, $t_n$.  For this
reason, a recurrence having the form~(\ref{eq3}) is called a
\emph{tree recurrence}.  Tree recurrences are the main object of study
of this paper.

One way to view the tree recurrence~(\ref{eq3}) is as a linear
transformation on sequences
\begin{equation}\label{eq4}
(f_n)= {\cal K}\left[(t_n)\right],
\end{equation}
that takes a toll sequence $(t_n)$ and returns the corresponding
average-cost sequence $(f_n)$. The functional $\cal K$ is fully
determined by the splitting probabilities $p_{n,k}$. A classical
approach to the derivation of explicit forms consists in introducing
\emph{generating functions} (GFs). Fix a sequence of
\emph{normalization constants} $\omega_n$ (that are problem-specific)
and define the generating functions
\[
f(z):=\sum f_n \omega_n z^n, \qquad t(z):=\sum t_n \omega_n z^n.
\]
Then, the transformation~$\cal K$ induces another linear
transformation $\cal L$ on GFs:
\begin{equation}\label{eq5}
f(z)= {\cal L}\left[t(z)\right].
\end{equation}
With an adequate choice of the constants~$\omega_n$, explicit forms
of~$f_n$ can often be obtained, provided at least that toll sequences
are of a simple enough form.

Our main objective is to develop generating-function methods
by which one can quantify the 
way the \emph{asymptotic form} of (expected) costs relates to
properties of the \emph{toll sequence}.
It is known that
asymptotic properties of number sequences (as the index~$n$ tends to
infinity) are closely related to the nature of
the~\emph{singularities} of the corresponding generating functions.
This suggests that we examine the way the operator~$\cal L$
operates on
scales of singular functions and view it as a ``\emph{singularity
  transformer}''.  Informally, there is a
transformation~$\widehat{\cal L}$, induced by~$\cal L$ and acting on
an asymptotic scale of functions singular at some fixed point~$z_0$.
Using $\operatorname{Sing}(f(z))$ to denote the expansion of~$f(z)$ at
the singularity~$z_0$, one has
\begin{equation}\label{eq6}
\operatorname{Sing}(f(z))=\widehat{\cal L}\left[t(z)\right]
\end{equation}
Under fairly general conditions, there is a tight coupling between
singular expansions of a generating function and the asymptotic form
of its coefficients.  The outcome of this process, justified by
\emph{singularity analysis}~\cite{FlOd90b,Odlyzko95}, is a direct
relation written figuratively as
\begin{equation}\label{eq7}
\operatorname{Asympt} ((f_n))=\widehat{\cal
K}\left[(t_n)\right],
\end{equation}
where~$\widehat{\cal K}$ depends on $(t_n)$  via the 
structure of its generating function~$t(z)$.

The path we follow in this paper is the one given
by~(\ref{eq4})--(\ref{eq7}), which is then globally
summarized by the following diagram:
\begin{equation}\label{commdiag}
\begin{array}{ccccc}
\ds (t_n) \mathop{\longrightarrow}^{\ds\cal K} (f_n) &&&&
\ds (t_n) \mathop{\longrightarrow}^{\ds\widehat{\cal K}}
\operatorname{Asympt}((f_n)) \\
\\
\ds \Downarrow &&&&\ds \Uparrow
\\
\\
\ds t(z) \mathop{\longrightarrow}^{\ds\cal L} f(z) && \Longrightarrow &&
\ds t(z) \mathop{\longrightarrow}^{\ds\widehat{\cal L}} 
 \operatorname{Sing}(f(z))
\end{array}
\end{equation}
We propose to develop a collection of generic tools that
supplement the basic singularity analysis framework of Flajolet and
Odlyzko~\cite{FlOd90b}.  In particular, we discuss in the next
sections the action on singularities of differential and integral
operators, as well as of Hadamard products.  As a result, the
way~$\cal L$ operators associated with many recurrences transform
singularities can be analyzed precisely. This in turn yields a fairly
general classification of the asymptotic growth phenomena associated
to a variety of classical tree recurrences, including the ones of
binary search trees, binary trees, and union--find trees, which will
serve here as guiding examples.

\medskip

Most of the existing computer science literature is devoted to the
``deterministic'' divide-and-conquer recurrences that correspond to a
splitting size $K_n$ that is deterministic, depending on~$n$
alone---typically, $K_n\equiv\lfloor n/2\rfloor$. In such a case, the
probability distribution $(p_{n,k})_{k=0}^n$ is supported at a single
point.  The main asymptotic order of~$f_n$ is then given by what
Cormen, Rivest, and Leiserson have termed ``master theorems'':\
see~\cite{CoLeRi90,GrKn82,SeFl96}.  (Usually, the finer 
characteristics of the asymptotic regime involve fractal
fluctuations~\cite{FlGo94,SeFl96}.) What we consider here instead are
methods for dealing with ``stochastic'' divide-and-conquer
recurrences, where $K_n$ is a random variable (dependent on~$n$) with
support spread over a whole subinterval in~$(0,n)$.  This stochastic
case is discussed by Roura in~\cite{Roura97}: Roura's arguments are
based on elementary real analysis, so that they are of quite a wide
scope, but his estimates are by nature mostly confined to first-order
asymptotics. In this article we show that, in the many cases of
practical interest where some strong complex-analytic structure is
present, full asymptotic expansions can be derived. Our treatment is
somewhat parallel in spirit to that of Knuth and Pittel whose
inspiring work~\cite{KnPi89} provided one of the initial
motivations\footnote{ See also Pittel's interesting recent
  article~\cite{Pittel99} which appeared while our own work was still
  in progress.  } for the present study. An additional benefit of the
complex-analytic approach is that it often gives access to variances
and higher moments, in which case  the limit distribution of
costs can be identified.

\section[Special tree recurrences]{Some ``special'' tree recurrences}\label{sec:some-special-tree}

In this section, we briefly review some tree recurrences that
are of special interest in combinatorial mathematics and analysis of
algorithms.

\subsection{The binary search tree recurrence}
One of the simplest model of random trees is defined as follows: To
determine a tree $T_n$ of size~$n\ge1$, take a root and append to it a
left subtree of size $k$ and a right subtree of size~$n-k$, where $k$
is uniformly distributed over the set of permissible
values~$\{0,1,\ldots n-1\}$; a tree of size~$0$ is the empty tree.  In
earlier notations, this process corresponds to
\begin{equation}\label{bst1}
p_{n,k}\equiv \Pr(K_n = k) :=\frac{1}{n},\qquad\hbox{for}\quad
k=0,1,\ldots n-1.
\end{equation}
As is well known, the model defined by~\eqref{bst1} corresponds to
random trees defined by either the binary search tree data structure
or the quicksort
algorithm~\cite{Knuth98a,Mahmoud92,Mahmoud00,SeFl96,ViFl90}.  The
corresponding tree recurrence~\eqref{eq3} is then
\begin{equation}\label{bst2}
f_n=t_n+\frac{2}{n}\sum_{k=0}^{n-1} f_k,
\end{equation}
with $f_0 := t_0$.

Ordinary GFs determined by the choice of coefficients $\omega_n=1$ for
all~$n$ are then
\[
f(z):=\sum_{n\ge0} f_n z^n, \qquad t(z):=\sum_{k\ge0} t_n z^n,
\]
and standard rules for the manipulation of GFs translate~(\ref{bst2})
into a linear integral equation
\[
f(z)=t(z)+2\int_0^z f(w) \, \frac{dw}{1-w}.
\]
Differentiation yields the ordinary differential equation
\[
f'(z)=t'(z)+\frac{2}{1-z}f(z),
\] 
which is then solved by the variation-of-constants method:
\begin{equation}\label{bst3}
f(z)=\cal L[t(z)], \qquad
\hbox{where}\quad \cal L[t(z)]:=(1-z)^{-2}\int_0^z
\left(\partial_w t(w)\right) (1-w)^2\, dw.
\end{equation}
In~(\ref{bst3}) we have assumed without loss of generality the initial
conditions $t_0=f_0=0$ (thanks to linearity and the fact that the
transform of~$t_n \equiv \delta_{n,0}$ is $n+1$). The
notation~$\partial_w$ 
borrowed from differential algebra is used to denote derivatives
whenever the operator nature of transformations is to be stressed.

It is instructive to follow what Greene and Knuth call the
``repertoire'' approach~\cite{GrKn82}.  This consists in building a
repertoire of the ($\cal K$ or $\cal L$) transforms of basic tolls,
then trying to determine the effect of a new toll by expressing it in
the basis of known tolls. What is convenient here is the class of
tolls
\[
t_n^\alpha:=\binom{n+\alpha}{\alpha} = \frac{(\alpha + 1)(\alpha + 2)
  \cdots (\alpha + n)}{n!}, \qquad \hbox{i.e.,}\quad
t^\alpha(z)=(1-z)^{-\alpha-1}.
\]
Then, by~(\ref{bst3}) one finds, for~$\alpha\not=1$,
\[\begin{array}{lll}
f^{\alpha}(z)&=&\ds \frac{\alpha+1}{\alpha-1}\left[(1-z)^{-\alpha-1}-(1-z)^{-2}
\right]\\
f^\alpha_n&=&\ds
\frac{\alpha+1}{\alpha-1}\left[\binom{n+\alpha}{\alpha}-\binom{n+1}{1}\right],
\end{array}
\]
while $\alpha=1$ leads to
\[
f^{1}(z)=\frac{2}{(1-z)^2}\log\frac{1}{1-z}, \qquad
f^{1}_n=2(n+1)\left(H_{n+1}-1\right),
\]
with $H_n : =1+\tfrac12+\cdots+\tfrac1n$ the $n$th harmonic number.

Stirling's formula implies asymptotically, for~$\alpha$ not a negative
integer,
\[
\binom{n+\alpha}{\alpha} \sim \frac{n^{\alpha}}{\Gamma(\alpha+1)},
\]
with $\Gamma$ the Euler gamma function. Then what goes on
is summarized  by the following table:
\begin{equation}\label{bst999}
\hbox{\begin{tabular}{lll}
\hline\hline
$\ds t_n=\binom{n+\alpha}{\alpha}$, $\alpha>1$ 
&&
$f_n\sim \frac{\alpha+1}{\alpha-1}\frac{n^{\alpha}}{\Gamma(\alpha+1)}$\\
$t_n=n+1$ 
&& 
$f_n\sim 2 n\log n$\\
$\ds t_n=\binom{n+\alpha}{\alpha}$, $0<\alpha<1$
&&
$f_n\sim \frac{1+\alpha}{1-\alpha}n$.\\
\hline\hline
\end{tabular}}
\end{equation}
The discontinuity in the asymptotic regime of~$f$ at~$\alpha=1$, where
a logarithm appears, is noticeable.
Also, the tolls in the scale satisfying $t_n\ll n$ are seen to induce
costs that all collapse to linear functions.

A full discussion of the binary search tree recurrence necessitates
determining the effect of toll functions like $\sqrt{n}$, $\log n$,
and
$1/n^2$, a task which is not entirely elementary. By the remarks of
the introduction, this involves determining the singularities of the
corresponding generating functions \emph{and}, in view
of~(\ref{bst3}), making explicit the way singular expansions get
composed under differentiation and integration
(Section~\ref{sec:sing-expans-diff}). This subject will then be taken
up again in Section~\ref{sec:binary-search-tree}; the particular case
of the toll~$t_n=\log n$ is of special importance and will be treated
in detail there.

\subsection{The uniform binary tree recurrence} 
This recurrence is of the form ($n \geq 0$, with the convention~$f_0
  := t_0$)
\begin{equation}\label{cat1}
f_n=t_n+\sum_{k=0}^{n-1} \frac{C_k C_{n-1-k}}{C_n}\left(f_k+f_{n-k}\right),
\qquad 
\hbox{with}\quad
C_n :=\frac{1}{n+1}\binom{2n}{n},
\end{equation}
a Catalan number.  It corresponds to the uniform model of binary
trees, where all the $C_n$ binary trees with~$n$ internal nodes are
taken with equal likelihood. Indeed, the number of trees of size~$n$
satisfies the recurrence
\begin{equation}\label{cat2}
C_n=\sum_{k=0}^{n-1} C_k C_{n-1-k}\quad (n \geq 1), \qquad
C_0=1,
\end{equation} 
as seen from a root decomposition. The
quantity~$p_{n,k}=C_kC_{n-1-k}/{C_n}$ is then the probability that a
tree of size~$n$ has left and right subtrees of respective sizes $k$
and $n-1-k$.

The GF of Catalan numbers satisfies a relation that is the image of
the recurrence~(\ref{cat2}), namely, $C(z)=1+zC(z)^2$, so that
\begin{equation}\label{cat3}
C(z)=\frac{1}{2z}\left(1-\sqrt{1-4z}\right).
\end{equation}
In order to solve~(\ref{cat1}) by generating functions, one should use
as
normalization constants the quantities~$\omega_n=C_n$, and introduce
\begin{equation}\label{catcost1}
t(z):=\sum_{n\ge0} t_n C_n z^n, \qquad
f(z):=\sum_{n\ge0} f_n C_n z^n.
\end{equation}
Then~(\ref{cat1}) translates into a linear algebraic equation,
\[
f(z)=t(z)+2zC(z)f(z),
\]
from which the form of the~$\cal L$ operator immediately results:
\begin{equation}\label{catcost2}
f(z)=\cal L[t(z)], \qquad
\hbox{where}\quad \cal L[t(z)]=\frac{1}{\sqrt{1-4z}} t(z).
\end{equation}
This form makes it possible to analyze directly only a restricted
collection of tolls, for instance, ones of the form
$t_n^r:=(n+1)n\cdots(n-r+2)$ (by differentiation), or
$t_n^{-r}=1/((n+2)(n+3)\cdots(n+r+1))$ (by integration). However,
tolls of such simple forms as $\sqrt{n}$, $H_n$, and~$\log n$, are
left out of the scale of the $t_n^{\pm r}$.

Define the \emph{Hadamard product} of two entire series or two
functions analytic at the origin, $a$ and~$b$, as their termwise
product,
\begin{equation}\label{hadam-def}
a(z)\odot b(z)=\sum_{n\ge0} a_n b_n z^n,
\qquad
\hbox{if}\quad
a(z)=\sum_{n\ge0} a_n z^n,\quad
b(z)=\sum_{n\ge0} b_n z^n.
\end{equation}
Then, from~(\ref{catcost1}) and~(\ref{catcost2}), the cost functional
is expressed by the
modified transformation (of~$\cal L$ type)
\begin{equation}\label{cat99}
f(z)=\frac{\tau(z)\odot C(z)}{\sqrt{1-4z}}, \quad
\hbox{where}\quad
f(z)=\sum_n f_n C_nz^n, \quad\tau(z):=\sum_n t_n z^n\,.
\end{equation}
This now relates the \emph{ordinary} generating function $\tau(z)$ of
the tolls and the \emph{normalized} generating function $f(z)$ of the
costs (with the $\omega_n=C_n$ normalization) via a Hadamard product.

Determining the way costs get transformed under this model then
necessitates a way to combine singular expansions under Hadamard
products. This is the central part of our article; see
Section~\ref{sec:hadam-prod-transf}, where a general theorem is
stated.  The ``critical'' value for tolls at which a discontinuity in
the induced costs manifests itself is now at $t_n=\sqrt{n}$,
and\footnote{%
  The notation $x=\Theta(y)$ expresses the inequalities $c_1y<x<c_2y$
  for some constants $c_1,c_2$ satisfying $0<c_1<c_2<+\infty$.  }
\begin{equation}\label{cat999}
\hbox{\begin{tabular}{lll}
\hline\hline
$t_n=n^{\alpha}$, $\alpha>1/2$ && $f_n=\Theta(n^{\alpha+\frac12})$\\
$t_n=n^{1/2}$ && $f_n=\Theta(n\log n)$\\
$t_n=n^{\alpha}$, $0 < \alpha < 1/2$ && $f_n=\Theta(n)$.\\
\hline\hline 
\end{tabular}%
}
\end{equation}
This phenomenon observed in~\cite[Prop.~2]{FlSt87} [of
which~(\ref{cat999}) above corrects a few misprints] neatly
distinguishes the binary Catalan model from the binary search tree
model, as seen by comparing~(\ref{cat999}) to~(\ref{bst999}).  A proof
accompanied by complete expansions will be given in the application
section: see Section~\ref{sec:binary-tree-recurr} below.

\subsection{The union--find tree recurrence} 
\label{sec:cayl-tree-recurr-1}
By a result attributed to Cayley, there are $U_n =n^{n-2}$ ``free''
unrooted trees (i.e., labeled connected acyclic graphs) on~$n$ nodes,
and, accordingly, $T_n=n^{n-1}$ rooted trees.  Consider the model in
which initially each unrooted tree of size~$n$ is taken with equal
likelihood.  Choose an edge at random amongst any of the possible
$n-1$ edges of the tree, orient it in a random way, then cut it.  This
separates the tree into an ordered pair of smaller trees that are now
rooted.  Continue the process with each of the resulting subtrees,
discarding the root.  Assume that the cost incurred by selecting the
edge and splitting the tree is $t_n$.  Then the total cost incurred
when starting from a random unrooted tree and recursively splitting it
till the completely disconnected graph is obtained satisfies the
recurrence $(n \geq 1)$
\begin{equation}\label{cay1}
f_n=t_n+\sum_{k=1}^{n-1}p_{n,k}(f_k+f_{n-k}), 
\quad\hbox{where}\quad
p_{n,k}=\binom{n}{k}\frac{k^{k-1}(n-k)^{n-k-1}}{2(n-1)\, n^{n-2}}.
\end{equation}
(Proof: There are $n^{n-1}$ rooted trees on~$n$ nodes and the binomial
coefficient takes care of relabellings.)  The recurrence~(\ref{cay1})
has been studied in great detail by Knuth and Pittel in~\cite{KnPi89},
an article that largely motivated our study.  In fact, there are good
algorithmic reasons for considering the recurrence~(\ref{cay1}): if
time is reversed, then the recursion describes the evolution of a
random graph from totally disconnected to tree-like, when successive
edges are added at random. The latter is exactly the probabilistic
model involved in the ``union--find'' (or equivalence-finding)
algorithm~\cite{CoLeRi90,Sedgewick88,ViFl90}, for which detailed
analyses had been provided by Knuth and
Sch\"onhage~\cite{KnSc78} in 1978\footnote{%
  Precisely, the model is known as the ``random spanning tree model''.
  The derivation of our equation~(\ref{cay2}) closely mimics
  Section~11 of~\cite{KnSc78}.}.  (Note that this model is not the same
as the simply generated family of Cayley trees.)

Let $T(z)$, $U(z)$ be the exponential generating functions of the
sequences~$(T_n)$, $(U_n)$, that is,
\[
T(z)=\sum_{n=1}^\infty n^{n-1}\frac{z^n}{n!}, \qquad
U(z)=\sum_{n=1}^\infty n^{n-2}\frac{z^n}{n!},
\]
It is a well-known fact of combinatorics that $T(z)$ satisfies the
functional relation $T(z)=ze^{T(z)}$, and one has $U=T-(T^2/2)$;
see~\cite{GoJa83,Knuth97,Stanley98}.  Define now the generating
functions
\[
t(z)=\sum_{n\ge1}t_n(n^{n-1}-n^{n-2})\frac{z^n}{n!},\qquad
f(z)=\sum_{n\ge1}f_n n^{n-1}\frac{z^n}{n!},
\]
where the normalization constants for~$f(z)$ are $\omega_n=n^{n-1}/n!$
and, for convenience, a marginally different normalization,
$\omega'_n=n^{n-2}(n-1)/n!$, has been introduced in the case
of $t(z)$.
Then the recurrence~(\ref{cay1}) has the form of a binomial
convolution, so that the cost GF~$f(z)$ satisfies
\[
f(z)-\int_0^z f(w)\, \frac{dw}{w}=t(z)+f(z)T(z).
\]
By differentiation, this last relation transforms into a linear
differential equation of the first order, itself readily
solved by the
variation-of-constants method. Assuming (without loss of generality)
the initial condition $t_1=f_1=0$, the solution found is
\begin{equation}\label{cay2}
f(z)=\frac{T(z)}{1-T(z)}\int_0^z \partial_w t(w)\, \frac{dw}{T(w)}.
\end{equation}

In terms of the ordinary generating function of costs, namely,
\[
\tau(z):=\sum_{n\ge2}t_n z^n,
\]
equation~(\ref{cay2}) can be rephrased as an integral transform
involving a Hadamard product, namely,
\begin{equation}\label{cay3}
f(z)=\frac12\frac{T(z)}{1-T(z)}\int_0^z \partial_w \left(\tau(w)\odot
T(w)^2\right) \, \frac{dw}{T(w)}.
\end{equation}
The dominant singularity at~$z=e^{-1}$ of the Cayley tree function~$T$
is
well known to be of the square root type. Then the integral
transform~(\ref{cay3}) operates in a way that combines a Hadamard
product and ordinary products, as well as integration and
differentiation.  This subject will be resumed in
Section~\ref{sec:cayl-tree-recurr}, after general theorems have been
established by which one can cope with such situations.  The final
conclusions turn out to be qualitatively similar to what was observed
for the Catalan model in~(\ref{cat999}).

\section{Singular expansions, differentiation, and integration}\label{sec:sing-expans-diff}

Singularities of generating functions encode very precise information
regarding the asymptotic behaviour of coefficients. In this section,
we first recall in Subsection~\ref{subsec:basicSA} the principles of a
process by which this information can be extracted: this is the
singularity analysis framework of~\cite{FlOd90b,Odlyzko95}. We then
prove that functions amenable to singularity analysis are closed under
integration and differentiation; see Subsection~\ref{subsec:diffint}.
These operations have already been seen to intervene in the
analysis of some of the major tree recurrences.

\subsection{Basics of singularity analysis}\label{subsec:basicSA}

Singularity analysis deals with functions that have isolated
singularities on the boundary of their disc of convergence and are
consequently continuable to  wider areas of the complex plane.  The
case of a unique dominant singularity suffices for the applications
treated here. (In addition, the case of finitely many dominant
singularities is easily reduced to this situation by using composite
contours and cumulating contributions arising from individual
singularities.) Given the obvious scaling rule,
\[
[z^n] f(z)=\rho^{-n}f(\rho z),
\]
one may restrict attention, whenever necessary, to the case where the
singularity is at~1. The scaling rule shows that the position of the
singularity [at~$\rho$ for $f(z)$] introduces an exponential scaling
factor ($\rho^{-n}$) multiplied by the coefficient of a function
singular at~1 [the function~$f(\rho z)$].

\begin{definition}\label{def:1}\em
    A function defined by a Taylor series with radius of convergence
  equal to~$1$ is \emph{$\Delta$-regular} if
  it can be analytically continued in a domain
\begin{equation*}
   \Delta(\phi,\eta) := \{z: |z| < 1 + \eta, |\Arg(z-1)| > \phi\},
\end{equation*}
for some $\eta > 0$ and $0 < \phi < \pi/2$. A function $f$ is said to
admit a \emph{singular expansion} at $z=1$ if it is $\Delta$-regular
and
\begin{equation}
   \label{eq:27}
   f(z) = \sum_{j=0}^J c_j(1-z)^{\alpha_j} + O(|1-z|^A)
\end{equation}
uniformly in $z \in \Delta(\phi,\eta)$,
for a sequence of complex numbers
$(c_j)_{0 \leq j \leq J}$ and an increasing sequence of real numbers
$(\alpha_j)_{0 \leq j \leq J}$ satisfying $\alpha_j < A$.  It is said
to satisfy a singular expansion \emph{``with logarithmic terms''} if,
similarly, 
\begin{equation}\label{eq:27ext}
f(z) = \sum_{j=0}^J c_j\left(L(z)\right)(1-z)^{\alpha_j} + O(|1-z|^A),
\qquad
L(z):=\log\frac{1}{1-z},
\end{equation}
where~each $c_j(\cdot)$ is a polynomial.
\end{definition}
Note that, by assumption, the $O(\cdot)$ error term in~\eqref{eq:27}
must hold uniformly in $z \in \Delta(\phi,\eta)$.
We also allow in the
usual way infinite asymptotic expansions
representing an infinite collection of mutually compatible expansions
of type~\eqref{eq:27}.

For the sake of notational simplicity, we shall mostly limit our
statements to the basic case~(\ref{eq:27}) and briefly comment on how
they extend to the logarithmic case~(\ref{eq:27ext}).  The basic
theorem is the following: 
\begin{theorem}[Basic singularity analysis~\cite{FlOd90b}]\label{thm:basic-sa} 
  If~$f(z)$ admits a singular expansion of the form~(\ref{eq:27})
  valid in a $\Delta$-domain, then 
\begin{equation}\label{xfer}
[z^n]f(z)=\sum_{j=0}^J
c_j\binom{n-\alpha_j-1}{-\alpha_j-1}+O(n^{-A-1}).
\end{equation}
\end{theorem}
(The proof of this and similar results is based on an extensive use of
Hankel contours; see the already cited references.)  The last
expansion can be rephrased as a standard asymptotic expansion since,
for~$\alpha\not\in\{0,1,2\ldots\}$, one has
\[
\binom{n-\alpha-1}{-\alpha-1}\sim
\frac{n^{-\alpha-1}}{\Gamma(-\alpha)}
\left(1+\frac{\alpha(\alpha+1)}{2n}+
  \frac{\alpha(\alpha+1)(\alpha+2)(3\alpha+1)}{24n^2}+\cdots\right),
\]
while all the terms corresponding to~$\alpha$ a nonnegative integer
have an asymptotically null contribution. When logarithmic terms are
present in the singular expansion, corresponding logarithmic terms
arise in the asymptotic expansion of coefficients. The calculations
are conveniently carried out by differentiation with respect to the
parameter~$\alpha$:
\[
[z^n](1-z)^{\alpha}L(z)^r=(-1)^r\frac{\partial^r}{\partial
  \alpha^r}[z^n](1-z)^\alpha =(-1)^r\frac{\partial^r}{\partial
  \alpha^r}\binom{n-\alpha-1}{-\alpha-1},
\]
which yields for instance ($\alpha\not\in\{0,1,2\ldots\}$):
\[
\begin{array}{lll}
\ds [z^n](1-z)^{\alpha}L(z)& =&\ds -\frac{\partial}{\partial \alpha}
\binom{n-\alpha-1}{-\alpha-1}
\\
&=&\ds \ds \binom{n-\alpha-1}{-\alpha-1}
\left(\frac{1}{-\alpha}+\frac{1}{1-\alpha}+\cdots+\frac{1}{n-1-\alpha}
\right)\\
&=& \ds\frac{n^{-\alpha-1}}{\Gamma(-\alpha)}\left(\log
n-\psi(-\alpha)+O\left(\frac{\log n}{n}\right)
\right).
\end{array}\]
(Here $\psi$ is the logarithmic derivative of~$\Gamma$.)

The same proof techniques also make it possible to translate error
terms involving logarithmic terms; see~\cite{FlOd90b} for details.  In
particular, the following transfer holds for~$A$ and~$B$ real numbers:
\begin{equation}\label{log-xfer}
O((1-z)^A L^B(z))\qquad\leadsto\qquad O(n^{-A-1}\log^B n).
\end{equation}

Finally, we shall make use of a result which renders amenable to
singularity analysis generating functions whose coefficients involve
powers of $n$ and its logarithms.

\begin{definition}\em
  The generalized polylogarithm~$\Li_{\alpha,r}$, where~$\alpha$ is an
  arbitrary complex number and $r$ a nonnegative integer is defined
  for~$|z|<1$ by 
\[
\Li_{\alpha,r}(z) :=\sum_{n\ge1} (\log n)^r \frac{z^n}{n^\alpha},
\]
and the notation~$\Li_\alpha$ abbreviates~$\Li_{\alpha,0}$.
\end{definition}
In particular, one has $\Li_{1,0}(z)=\Li_1(z)=L(z)$, the usual
logarithm, cf.~(\ref{eq:27ext}).  The singular expansion of the
polylogarithm, taken from~\cite{Flajolet99}, involves the Riemann
$\zeta$ function:
\begin{theorem}[Singularities of polylogarithms~\cite{Flajolet99}]\label{thm:li-sing}
  The function~$\Li_{\alpha,r}(z)$ is $\Delta$--continuable and, for
  $\alpha\not\in\{1,2,\ldots\}$, it satisfies the singular expansion
\begin{equation}\label{li-sing}
\Li_{\alpha,0}(z)\sim\Gamma(1-\alpha)w^{\alpha-1}+
\sum_{j\ge0}\frac{(-1)^j}{j!}\zeta(\alpha-j)w^j, \quad
w=-\log z=\sum_{\ell=1}^\infty \frac{(1-z)^\ell}{\ell}.
\end{equation}
For~$r>0$, the singular expansion of~$\Li_{\alpha,r}$ is obtained by
\[
\Li_{\alpha,r}(z)=(-1)^r\frac{\partial^r}{\partial
  \alpha^r}\Li_{\alpha,0}(z),
\]
and corresponding termwise differentiation of~\eqref{li-sing} with
respect to~$\alpha$.
\end{theorem}
In particular, for~$\alpha<1$, the main asymptotic term
of~$\Li_{\alpha,r}$ is
\[
\Gamma(1-\alpha) (1-z)^{\alpha-1} L^r(z).
\] 
Similar expansions hold when $\alpha$ is a positive integer;
see~\cite{Flajolet99} for details.

\begin{example} \emph{Stirling's formul\ae.}
  The factorial function, is attainable via the form
\[
\log n! =\log 1+\log 2+\cdots+\log n = [z^n]
\frac{1}{1-z}\Li_{0,1}(z),
\]
to which singularity analysis can be applied now that we have taken ordinary
generating functions. Theorem~\ref{thm:li-sing} yields
the singular expansion 
\[
\frac{1}{1-z}\Li_{0,1}(z)\sim\frac{L(z)-\gamma}{(1-z)^2}
+\frac12\frac{-L(z)+\gamma-1+\log{2\pi}}{1-z}+\cdots,
\]
from which Stirling's formula can be read off, by
Theorem~\ref{thm:basic-sa}: 
\[
\log n!\sim n\log n -n +\frac{1}{2}\log n+\log\sqrt{2\pi}+\cdots\,.
\]
[Stirling's constant~$\log\sqrt{2\pi}$ comes out as $-\zeta'(0)$.]
Similarly, the ``superfactorial function'',
\[
S(n) := 1^1\cdot 2^2\cdots
n^n\equiv\frac{(n!)^{n+1}}{1!\,2!\,\cdots\,n!},
\]
satisfies
\[
\log S(n) = [z^n]\frac{1}{1-z}\Li_{-1,1}(z),
\]
which gives rise to a second-order ``Stirling's formula'',
\[
S(n) \sim n^{\frac12n^2+\frac12n+\frac{1}{12}}e^{-\frac14 n^2} A,
\]
with
\[
A:= \exp \left(\frac1{12} - \zeta'(-1)\right) =
\exp\left(-\frac{\zeta'(2)}{2\pi^2}+\frac{\log(2\pi)+\gamma}{12}\right).
\]
(This last expansion, originally due to Glaisher, Jeffery, and Kinkelin,
goes back to the 1860s and it can be established by
Euler--Maclaurin summation; see Finch's book~\cite{Finch03} for
context and references.) The systematic character of the derivation
given here clearly applies to many similar functions.\hfill~$\Box$
\end{example}

Methods of the last example may be used more generally to determine
the Euler--Maclaurin constant relative to sums of the form $\sum (\log
n)^r/n^s$. The derivation by singularity analysis is quite systematic
and several formul{\ae} of Ramanujan can be obtained in this way, for
instance,
\[
\lim_{N\to\infty} \left(\sum_{n=1}^N \frac{\log^k
    n}{n}-\frac{\log^{k+1}N}{k+1}\right)=A_k, ~\hbox{with}~
A_k:=\frac{(-1)^k}{k+1}\frac{d^{k+1}}{ds^{k+1}}\left((s-1)\zeta(s)\right)_{s=1}
,\] involving the Stieltjes constants~$A_k$. See Berndt's account of
the problem in~\cite[p.~164]{Berndt85} and references therein.

%
%

\subsection{Differentiation and integration}\label{subsec:diffint}

In preparation for our later treatment of Hadamard products, we need a
theorem that enables us to differentiate local expansions of analytic
functions around a singularity. Such a result cannot of course be
unconditionally true; see, for example,~\eqref{eq:56}. However, it
turns out that functions amenable to singularity analysis
satisfy this property. The statement that follows is an adaptation
suited to our needs of well-known differentiability properties of
complex asymptotic expansions (see especially Theorem~I.4.2 of Olver's
book \cite[p.~9]{Olver74}).
\begin{theorem}[Singular differentiation]
   \label{thm:2A}
   If $f(z)$ is $\Delta$-regular and admits a singular expansion near
   its singularity in the sense of~\eqref{eq:27}, then for each
   integer $r > 0$, $\tfrac{d^r}{dz^r}f(z)$ is also $\Delta$-regular
   and admits an expansion obtained through term-by-term
   differentiation:
   \begin{equation*}
     \frac{d^r}{dz^r} f(z) = (-1)^r\sum_{j=0}^J c_j 
     \frac{\Gamma(\alpha_j+1)}{\Gamma(\alpha_j+1-r)} 
(1-z)^{\alpha_j-r} + O( |1-z|^{A-r}).
   \end{equation*}
\end{theorem}
\begin{proof}
Clearly, all that is required is to establish the effect of
  differentiation on error terms, which is expressed symbolically as
   \begin{equation*}
     \frac{d}{dz} O( |1-z|^A ) = O( |1-z|^{A-1} ).
   \end{equation*}
   By iteration, only the case of a single differentiation ($r=1$)
   needs to be considered.

\begin{figure}
  \centering
  \includegraphics[width=5truecm]{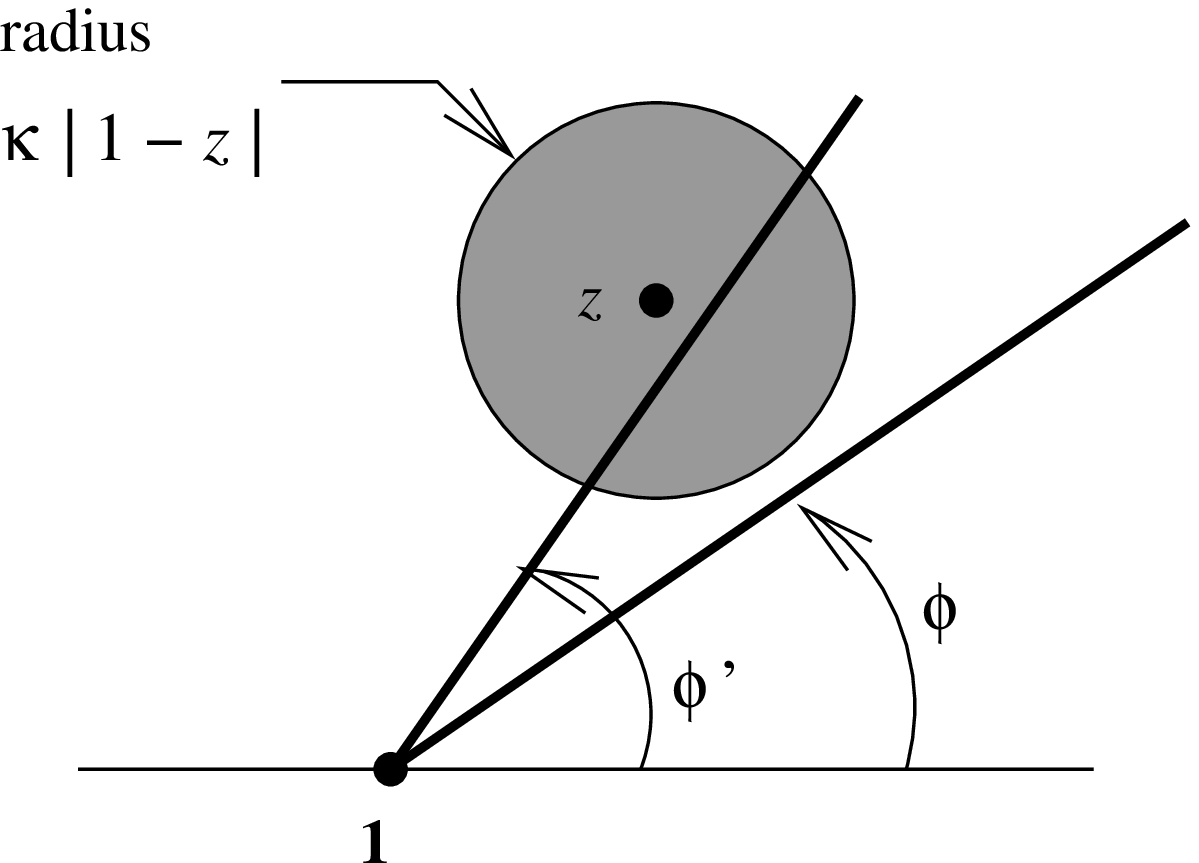} \hspace*{1truecm}
  \includegraphics[width=5truecm]{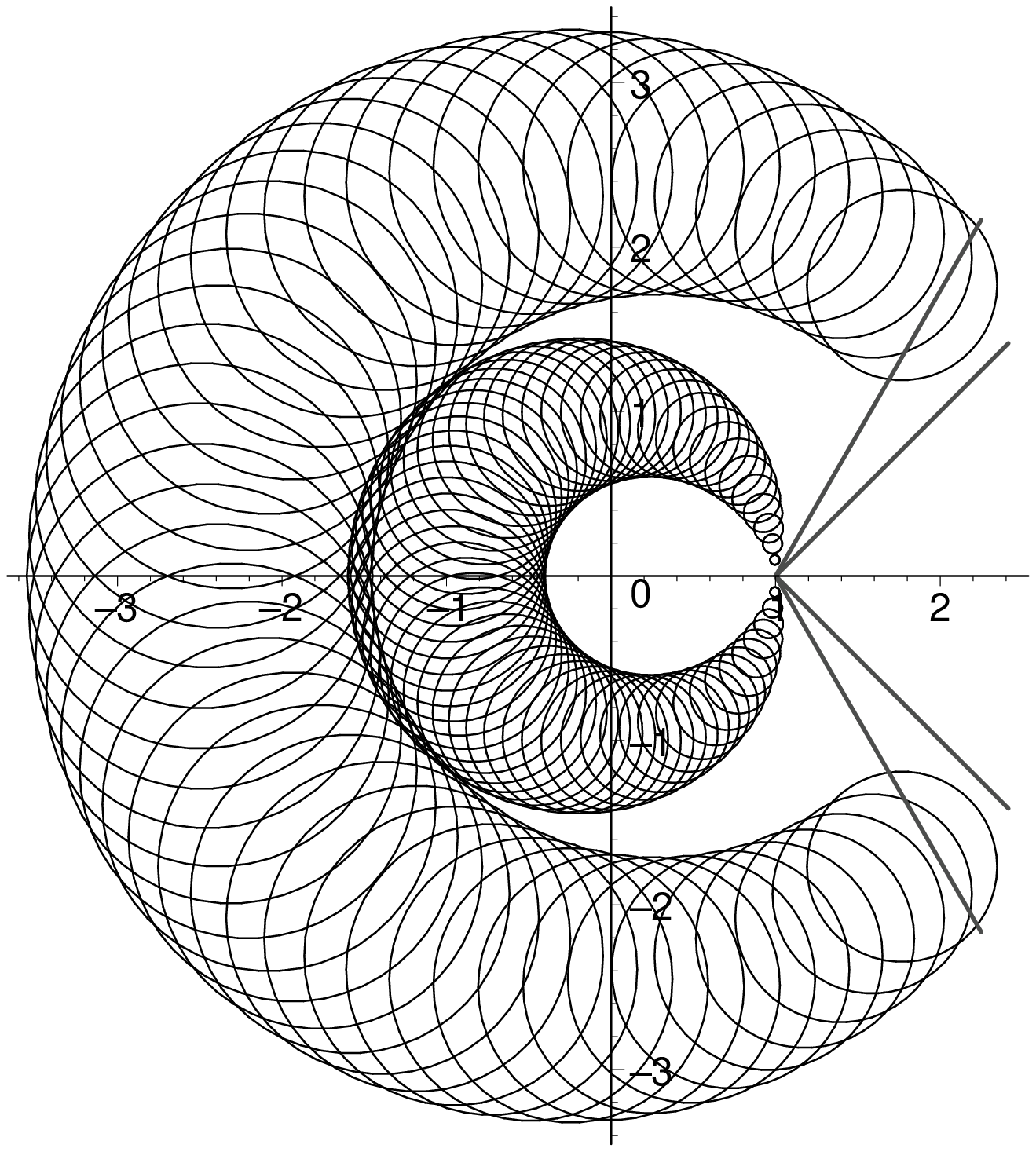}
      \caption{The contour $\gamma(z)$ used in the proof of the
        differentiation theorem: (left) the basic geometry;
(right) two sets  of circles $\gamma(z)$.}%
    \label{fig:diff}%
    \end{figure}
    
    Let~$g(z)$ be a function that is regular in a domain
    $\Delta(\phi,\eta)$ where it is assumed to satisfy
    $g(z)=O(|1-z|^A)$ for $z \in \Delta$.
    Choose a subdomain
$\Delta':=\Delta(\phi',\eta')$,
where $\phi<\phi'<\frac\pi2$ and $0<\eta'<\eta$. By elementary
geometry, for any sufficiently small~$\kappa>0$, the disc of radius 
$\kappa|z-1|$ centered at a value~$z\in\Delta'$
lies entirely in~$\Delta$; see Figure~\ref{fig:diff}.
 We fix such a small
value~$\kappa$ and let~$\gamma(z)$ represent the boundary
of that disc oriented positively.

The starting point is Cauchy's integral formula 
   \begin{equation}
     \label{eq:35}
     g'(z) = \frac{1}{2\pi i} \int_C g(w) \frac{dw}{(w-z)^2},
   \end{equation}
   a direct consequence of the residue theorem. Here~$C$  should
   encircle~$z$ while lying inside the domain of regularity of~$g$,
   and we opt for the choice $C\equiv\gamma(z)$. Then trivial bounds
   applied to~(\ref{eq:35}) give:
\[ \begin{array}{lll}
|g'(z)|&=& O\left(  |\!| \gamma(z)|\!|\cdot |1-z|^A
|1-z|^{-2}\right)
\\
&=& \ds O\left(|1-z|^{A-1}\right).
\end{array}
\]
The estimate involves the length of the contour, $|\!| \gamma(z)|\!|$,
which is $O(|1-z|)$ by construction, as well as the bound on~$g$
itself, which is $O(|1-z|^A)$ since all points of the contour are
themselves at a distance exactly of the order of $|1-z|$ from~1.
\end{proof}

For instance, taking
\[
g(z) = \cos{\log{\left(\frac{1}{1-z}\right)}} \quad \text{and} \quad
g'(z) = - \frac{1}{1-z}\sin{\log{\left(\frac{1}{1-z}\right)}},
\]
we correctly predict that $g(z) = O(1) \Rightarrow g'(z) = O(
|1-z|^{-1})$. On the other hand, the apparent paradox given by the
pair
\begin{equation}
  \label{eq:56}
g(z) = \cos{{\left(\frac{1}{1-z}\right)}} \qquad \text{and}\qquad
g'(z) = - \frac{1}{(1-z)^2}
\sin{{\left(\frac{1}{1-z}\right)}} ,
\end{equation}
is resolved by observing that in no nondegenerate sector around $z=1$
do we have $g(z) = O(1)$.

It is also well known that integration of asymptotic expansions is
usually easier than differentiation. Here is a statement
custom-tailored to our needs.

\begin{theorem}[Singular integration]  \label{thm:2B}
  Let $f(z)$ be $\Delta$-regular and admit a $\Delta$-expansion near
  its singularity in the sense  of~\eqref{eq:27}.  Then
  $\int_0^zf(t)\,dt$ is also $\Delta$-regular.  Assume that none of
  the quantities $\alpha_j$ and $A$ equals~$-1$.
  
  $(i)$~If~$A<-1$, then the singular expansion of $\int f$ is
\begin{equation}\label{eq:in1}
\int_0^z f(t)\,dt=-\sum_{j=0}^J
\frac{c_j}{\alpha_j+1}(1-z)^{\alpha_j+1}+O\left(|1-z|^{A+1}\right).
\end{equation}

$(ii)$~If~$A>-1$, then the singular expansion of $\int f$ is 
\[
\int_0^z f(t)\, dt=-\sum_{j=0}^J
  \frac{c_j}{\alpha_j+1}(1-z)^{\alpha_j+1} +L_0
+O\left(|1-z|^{A+1}\right),
\]
where the ``integration constant'' $L_0$ has the value
\[
L_0:=\sum_{\alpha_j<-1}\frac{c_j}{\alpha_j+1}+
\int_0^1\Big[f(t)-\sum_{\alpha_j<-1}c_j(1-t)^{\alpha_j}\Big]\, dt.
\]
\end{theorem}
\begin{remark*}
The case where either some~$\alpha_j$ or~$A$ is $-1$ is easily treated
by the additional rules 
\[
\int_0^z (1-t)^{-1}\, dt=L(z), \qquad \int_0^z O(|1-t|^{-1})\, dt =
O(L(z)).
\]
Similar rules consistent with elementary integration are applicable
for powers of logarithms: they are derived from the easy identities
(for $\alpha \ne -1$) 
\[
\int_0^z (1-t)^\alpha L^r(t)\, dt=(-1)^r \frac{\partial^r}{\partial
  \alpha^r} \int_0^z (1-t)^{\alpha}\,dt = (-1)^{r+1}
\frac{\partial^r}{\partial \alpha^r}
\frac{(1-z)^{\alpha+1}}{\alpha+1},
\]
for~$r$ a positive integer. Furthermore, the corresponding
$O$--transfers hold true.  (The proofs are simple modifications of the
one given below for the basic case.)
\end{remark*}
\begin{proof} The basic technique consists in integrating, term by
  term, the singular expansion of~$f$. We let~$r(z)$ be the remainder
  term in the expansion of~$f$, that is,
\[
r(z):=f(z)-\sum_{j=0}^J c_j(1-z)^{\alpha_j}.\] By assumption,
throughout the $\Delta$-domain one has,  for some positive
constant~$K$,
\[
|r(z)|\le K|1-z|^A.
\]

  \begin{figure}  
    \centering
    \includegraphics[width=6truecm]{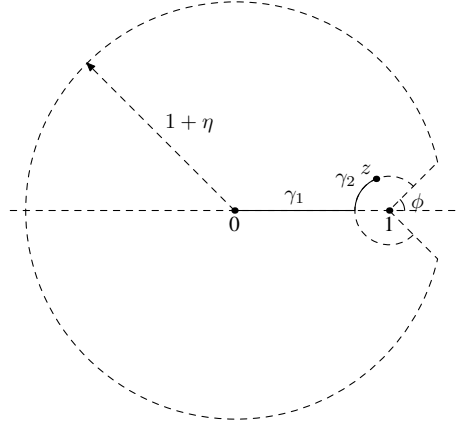}%
    \caption{The contour used in the proof of the integration theorem.}
    \label{fig:int}%
  \end{figure}
  
  $(i)$~\emph{Case $A<-1$.}  By straight-line integration between~$0$
  and~$z$, one finds~(\ref{eq:in1}), as soon as it has been
  established that
\[
\int_0^z r(t)\, dt=O\left(|1-z|^{A+1}\right).
\]
By Cauchy's integral formula, we can choose any path of integration
that stays within the region of analyticity of $r$. We choose the
contour $\gamma := \gamma_1 \cup \gamma_2$, shown in
Figure~\ref{fig:int}.
  Then\footnote{%
    The symbol $|dt|$ designates the differential line element~(often
    denoted by $ds$) in the corresponding curvilinear integral.  }
  \begin{align*}
    \left|\int_{\gamma} r(t) \,dt \right|
    &\leq \left| \int_{\gamma_1} r(t) \,dt \right| + \left|
      \int_{\gamma_2} r(t) \,dt \right|\\
    &\leq K \int_{\gamma_1} |1-t|^A \,|dt| + K
    \int_{\gamma_2} |1-t|^A|\,|dt|\\
    & = O(|1-z|^{A+1}).\\
  \end{align*}
  Both integrals are~$O(|1-z|^{A+1})$: for the integral
  along~$\gamma_1$, this results from explicitly carrying out the
  integration; for the integral along $\gamma_2$, this results from
  the trivial bound~$O(|\!|\gamma_2|\!| (1-z)^A)$.

  $(ii)$~\emph{Case $A>-1$.} We let $f_-(z)$ represent the 
  ``divergence part'' of~$f$ that gives rise to nonintegrability:
\[
f_-(z):=\sum_{\alpha_j<-1}c_j(1-z)^{\alpha_j}.
\]
Then with the decomposition $f=[f-f_-]+f_-$, integrations can be
performed separately.  First, one finds 
\[
\int_0^z f_-(t)\,
dt=-\sum_{\alpha_j<-1}\frac{c_j}{\alpha_j+1}(1-z)^{\alpha_j+1}
+\sum_{\alpha_j<-1}\frac{c_j}{\alpha_j+1}.
\]
Next, observe that the asymptotic condition guarantees the existence
of $\int_0^1$ applied to $[f-f_-]$, so that
\[
\int_0^z \left[f(t)-f_-(t)\right]\,dt =\int_0^1
\left[f(t)-f_-(t)\right]\,dt +\int_1^z \left[f(t)-f_-(t)\right]\,dt.
\]
The first of these two  integrals is a constant that
contributes to~$L_0$.
As to the second integral, term-by-term integration yields 
\[
\int_1^z \left[f(t)-f_-(t)\right]\,dt=-\sum_{\alpha_j>-1}
\frac{c_j}{\alpha_j+1} (1-z)^{\alpha_j + 1}+\int_1^z r(t)\, dt.
\]
The remainder integral is finite, given the growth condition on the
remainder term, and, upon carrying out the integration along the
rectilinear segment joining~$1$ to~$z$, trivial bounds show that it is
indeed~$O(|1-z|^{A+1})$.
\end{proof}

%
%

\section{Hadamard products and transformation of
  singularities}
\label{sec:hadam-prod-transf}
In this section we propose to examine the way singular expansions get
composed under Hadamard products defined 
at~(\ref{hadam-def}).  The Hadamard product is a bilinear form.
So if we have a set of functions admitting known singular expansions,
we need to establish their composition law, and this will give
composition rules for finite terminating expansions
(Subsection~\ref{subsec:had-pow}).  In order to extend this to
asymptotic expansions with error terms, we need to establish a theorem
providing the shape of
\[
O( |1-z|^A ) \odot O( |1-z|^B ).
\]
This is the more demanding part of the analysis, which is the subject
of Subsection~\ref{subsec:had-bigo}. Finally, in
Subsection~\ref{subsec:had-all}, we provide a summary
statement, Theorem~\ref{thm:hadmain}, to the effect that the class of
functions amenable 
to singularity analysis is  closed under Hadamard products and 
that the composition of singular expansions is effectively computable.

\subsection{Composition of singular elements}\label{subsec:had-pow}

The composition rule for polylogarithms is trivial, since
\[
\Li_{\alpha,r}(z) \odot \Li_{\beta,s}(z) = \Li_{\alpha+\beta,r+s}(z).
\]
However, polylogarithms do not have a simple composition rule with
respect to ordinary products.
We next turn to the composition rule for the basis formed by functions
of the form $(1-z)^{a}$, where $a$ may be any
real
number.  From the expansion
\begin{equation}
   \label{eq:39}
   (1-z)^{a} = 1 + \frac{-a}{1} z + 
\frac{(-a)(-a+1)}{2!} z^2 + \cdots
\end{equation}
around the origin,
we get through term-by-term multiplication
\begin{equation}
   \label{eq:40}
   (1-z)^{a} \odot (1-z)^{b} = \FF[-a,-b;1;z].
\end{equation}
Here $\FF$ represents 
the  classical \emph{hypergeometric function} of Gauss  defined 
by
\begin{equation}
   \label{eq:43}
   \FF[\alpha,\beta;\gamma;z] = 1 + \frac{\alpha\beta}{\gamma} \frac{z}{1!} +
   \frac{\alpha(\alpha+1)\beta(\beta+1)}{\gamma(\gamma+1)} 
\frac{z^2}{2!} + \cdots.
\end{equation}
From the transformation theory of hypergeometrics, see
e.g.~\cite[p.~163]{MR56:12235}, we
know that, in general, hypergeometric functions can
be expanded in the vicinity of $z=1$ by means of the $z \mapsto  1 - z$
transformation. Instantiation
of this transformation with $\gamma=1$ yields
\begin{multline}
   \label{eq:42}
   \FF[\alpha,\beta;1;z] = 
\frac{\Gamma(1-\alpha-\beta)}{\Gamma(1-\alpha)\Gamma(1-\beta)}
   \,\FF[\alpha,\beta;\alpha+\beta;1-z]\\
   + \frac{\Gamma(\alpha+\beta-1)}{\Gamma(\alpha)\Gamma(\beta)} 
(1-z)^{-\alpha-\beta+1}
   \,\FF[1-\alpha,1-\beta;2-\alpha-\beta;1-z].
\end{multline}
In other words, we can state the following proposition:
\begin{proposition}   \label{thm:3}
   When $a$, $b$, and $a + b$ are not integers, 
   the Hadamard product
   \[ (1-z)^{a} \odot (1-z)^{b} \]
   has an infinite $\Delta$-expansion with exponent scale
   \[
   \{0,1,2,\ldots\} \cup \{a + b + 1, a + b + 2,\ldots\},
   \]
   namely,
   \[
   (1-z)^{a} \odot (1-z)^{b} \sim \sum_{k \geq 0} \lambda_k^{(a,b)}
   \frac{(1-z)^k}{k!} + \sum_{k \geq 0} \mu_k^{(a,b)}
   \frac{(1-z)^{a + b +1+k}}{k!},
   \]
   where the
   coefficients
   $\lambda$ and $\mu$ are given by
   \begin{align*}
     \lambda_k^{(a,b)} &= \frac{\Gamma(1 + a + b)}
     {\Gamma(1+a)\Gamma(1+b)}
     \frac{\rising{(-a)}{k}\rising{(-b)}{k}}{\rising{(-a-b)}{k}},\\
     \mu_k^{(a,b)} &=
     \frac{\Gamma(-a-b-1)}{\Gamma(-a)\Gamma(-b)}
     \frac{\rising{(1+a)}{k}
       \rising{(1+b)}{k}}{\rising{(2+a+b)}{k}}.
   \end{align*}
   Here $\rising{x}{k}$ is defined when $k$ is a nonnegative
   integer as $x(x+1)\cdots(x+k-1)$.
\end{proposition}
\begin{remark*}
  \label{rem:a-b-integer}
  The case where either $a$ or $b$ is an integer  poses
  no difficulty: one has
\begin{itemize}
\item[---] $(1-z)^a \odot g(z)$ is a polynomial if $a=m$, where
  $m\in\Z_{\ge0}$;
\item[---] $(1-z)^a \odot g(z)$ is a derivative if $a=-m$ where
  $m\in\Z_{>0}$, since
\[
(1-z)^{-m}\odot g(z)=\frac{1}{(m-1)!}\partial_z^{m-1}
\left(z^{m-1}g(z)\right),
\]
and this case is covered by singular differentiation,
Theorem~\ref{thm:2A}.  
\end{itemize}
Notice that Proposition~\ref{thm:3} remains  valid in these two cases
with the natural convention that
$1/\Gamma(-j) = 0$ when $j \in \Z_{\geq 0}$. 

The case where $a + b \in \Z$ needs transformation formul{\ae}
that extend~(\ref{eq:42}) and are found explicitly in the books by
Abramowitz and Stegun~\cite[pp.~559--560]{AbSt73} and by Whittaker and
Watson~\cite[\S14.53]{WhWa27}.
\end{remark*}
\begin{remark*} 
  \label{rem:logarthmic-factors}
The case of expansions with logarithmic terms
 is covered
by ``differentiation under the integral sign'', 
as we now explain. Consider for instance
the Hadamard product 
\[
\left[(1-z)^{-\alpha}L(z)\right]\odot (1-z)^{-\beta}
= \frac{\partial}{\partial_\alpha}\, 
\FF[\alpha,\beta;1;z],
\]
where we assume for convenience that none of
$\alpha$, $\beta$, $\alpha + \beta$ is an integer. For any
fixed~$\beta$ and any fixed~$z$,
with, say, $z\in(0,1)$, both sides  of~(\ref{eq:42})
represent analytic functions  
of~$\alpha$. Thus, their derivatives 
with respect to~$\alpha$ are identical as functions 
of~$\alpha$. This induces a transformation formula, originally valid in the
stated $z$-range, which involves modified 
hypergeometric functions (these have additional $\psi$-factors
in their coefficients) obtained from the 
fundamental ${}_2F_1$ function by differentiation with respect to some
of the parameters. The modified  functions then do exist in extended 
regions of the complex $z$--plane as shown by taking the 
classical Barnes 
representations in terms of  contour integrals
(see, e.g.,~\cite[\S14.5]{WhWa27}) and then differentiating under the
integral sign. The net effect of this discussion is
that the fundamental transformation~(\ref{eq:42}) supports
differentiation with respect to~$\alpha,\beta$ and that the 
formally derived transformations provide analytically valid
composition formul{\ae} for Hadamard products 
\begin{equation}\label{hadlog}
[(1-z)^{-\alpha}L^k(z)] \odot [(1-z)^{-\beta} L^\ell(z)]
\end{equation}
 of the base functions.
\end{remark*}

In practice, for all the cases described above, one may often proceed
as follows: \textit{(i)}~take advantage of the a priori \emph{existence} of a
singular expansion of~$f\odot g$, with $f(z)=(1-z)^{a}$, 
$g(z)=(1-z)^{b}$ or some of their derivatives, that is valid
for~$z$ in a $\Delta$-region (here the slit complex plane);
\textit{(ii)}~compute an asymptotic expansion of the coefficients of $f\odot
g$ by multiplication of the asymptotic expansions of~$f_n$ and~$g_n$
as obtained via singularity analysis; \textit{(iii)}~reconstruct a singular
function that matches asymptotically $f_ng_n$ by using singularity
analysis in the reverse direction.  In 
Subsection~\ref{subsec:had-all}, this process is formalized by the
``Zigzag Algorithm''  and illustrated by the return of P\'olya's
drunkard.

\smallskip

Globally, we are facing a situation where polylogarithms are simple for
Hadamard products and relatively complicated for ordinary products,
with the dual situation occurring in the case of power functions.
Each particular situation is likely to dictate whether calculations
are best expressed in a basis of standard singular functions like
$\{(1-z)^{a}L(z)^k\}$ or  
with polylogarithms, $\{\Li_{\alpha,k}(z)\}$.

\subsection{Composition of error terms}\label{subsec:had-bigo}
We now examine how $O(\cdot)$ terms get composed
under Hadamard products. The task is easier when
the resulting function gets large at its singularity
as shown by Proposition~\ref{lem:3}. Fortunately,
thanks to the results of
Section~\ref{sec:sing-expans-diff}
regarding differentiation and integration, all cases can be reduced to this
one: see Proposition~\ref{thm:4} below.

{\sloppy
The starting point is  a general integral formula
due to Hadamard for \mbox{$(f \odot g)(z)$}, where 
\[
f(z) = \sum_{n \geq 0} f_nz^n \qquad \text{and} \qquad g(z) = \sum_{n \geq
   0} g_nz^n.
\]}%
Assume that $f$ and $g$ are analytic in the unit disc
and let $z$ be a complex number satisfying $|z| < 1$.
   Consider the integral
   \begin{equation}
     \label{eq:45}
     I = \frac{1}{2\pi i} \int_{\gamma_0} f(w)g\left(\frac{z}{w}\right)
     \frac{dw}{w},
   \end{equation}
taken (counterclockwise) along a contour~$\gamma_0$ 
which
is simply a circle of radius $\rho$ centered at the origin such that 
$|z| < \rho < 1$. In this way, both factors in the integrand are analytic
functions of~$w$ along the contour.
Evaluating the integral~(\ref{eq:45}) 
by expanding the functions, we find
\[
I = \sum_{n \geq 0} f_ng_n z^n.
\]
This is the classical formula of Hadamard for Hadamard products,
\begin{equation}
     \label{eq:45bis}
     (f \odot g)(z) = \frac{1}{2\pi i} \int_{\cal C}
     f(w)g\left(\frac{z}{w}\right)
     \frac{dw}{w},
   \end{equation}
valid, by analyticity, for any simple contour~$\cal C$ such that 
each~$w\in\cal C$ satisfies $|z|< |w| <1$.

\begin{proposition}   \label{lem:3}
   Assume that $f(z)$ and $g(z)$ are $\Delta$-regular in $\Delta(\psi_0, \eta)$
   and that 
   \[
   f(z) = O(|1-z|^a)\mbox{\rm \ \ and\ \ }g(z) = O(|1-z|^b), \quad
   z\in\Delta(\psi_0,\eta), 
   \]
   where $a$ and $b$
   satisfy $a+b+1 < 0$. Then the Hadamard product $(f \odot g)(z)$ is
   regular in a (possibly smaller)
$\Delta$-domain, call it~$\Delta'$, where it admits the expansion
   \begin{equation}
     \label{eq:44}
     (f \odot g)(z) = O( |1-z|^{a+b+1} ).
   \end{equation}
\end{proposition}



\begin{figure}
\begin{center}
\includegraphics[width=6truecm]{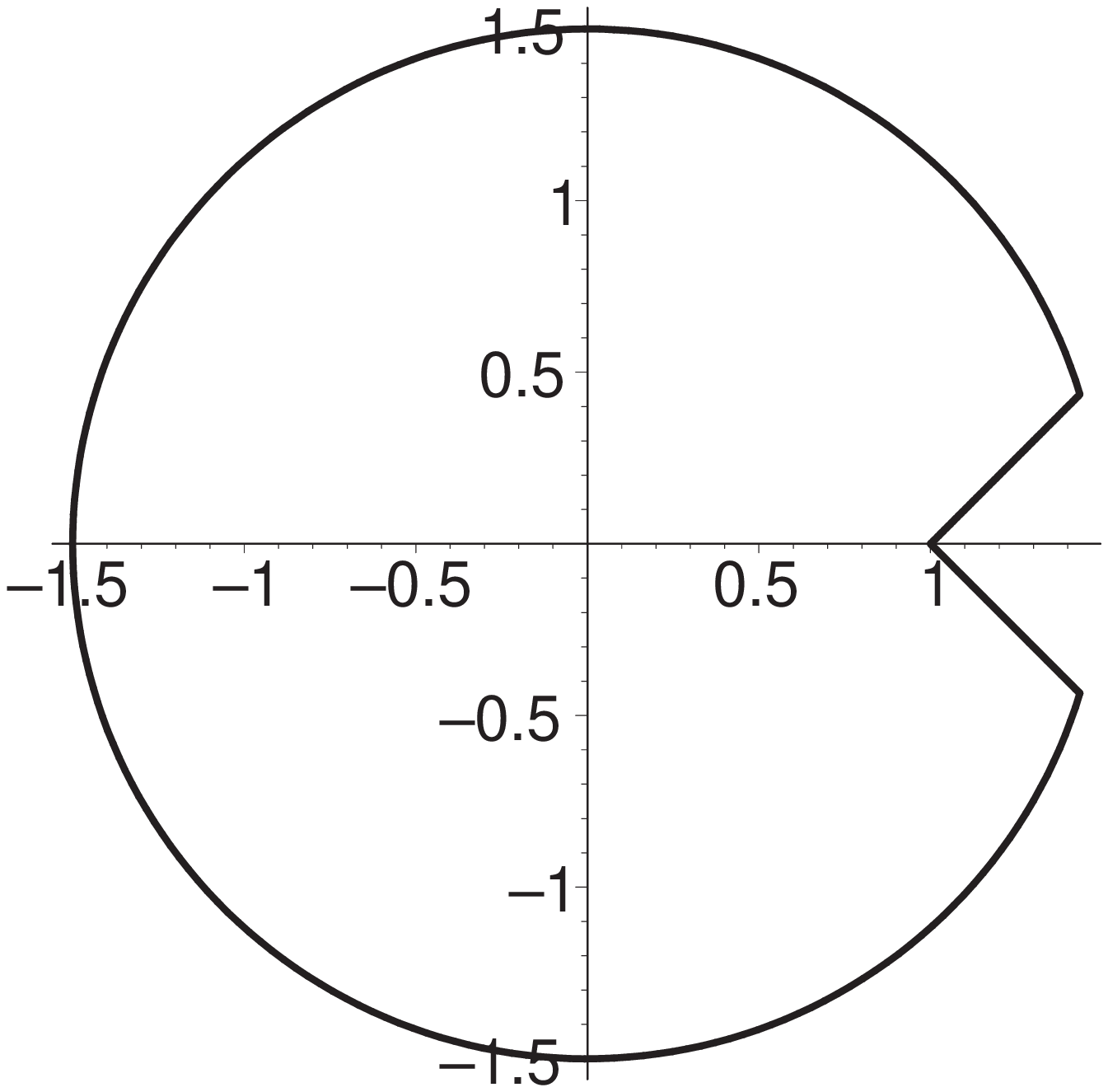}\hspace*{-1.5truecm}%
\raisebox{16pt}{\includegraphics[width=4.0truecm]{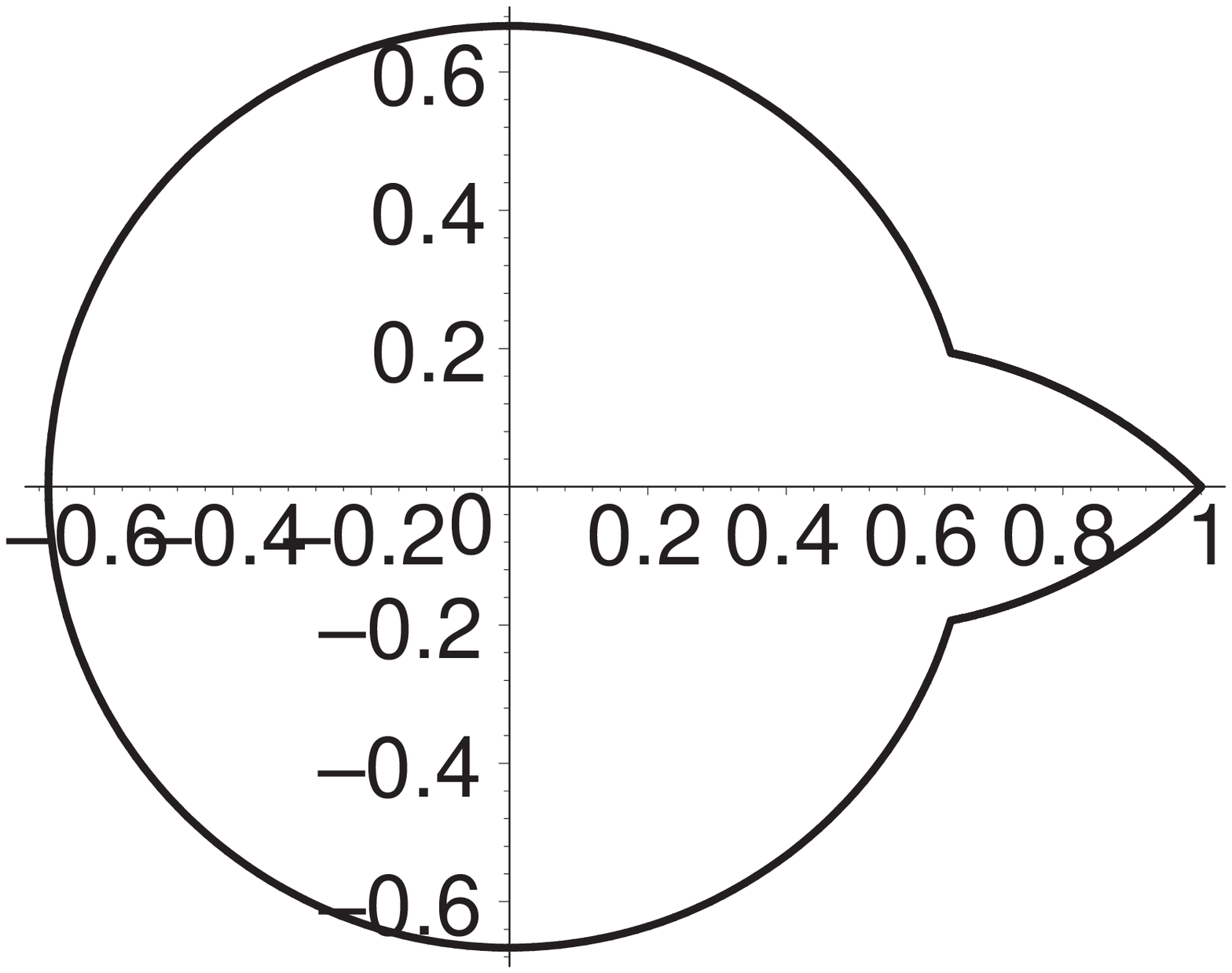}}\hspace*{-1.15truecm}%
\includegraphics[width=6truecm]{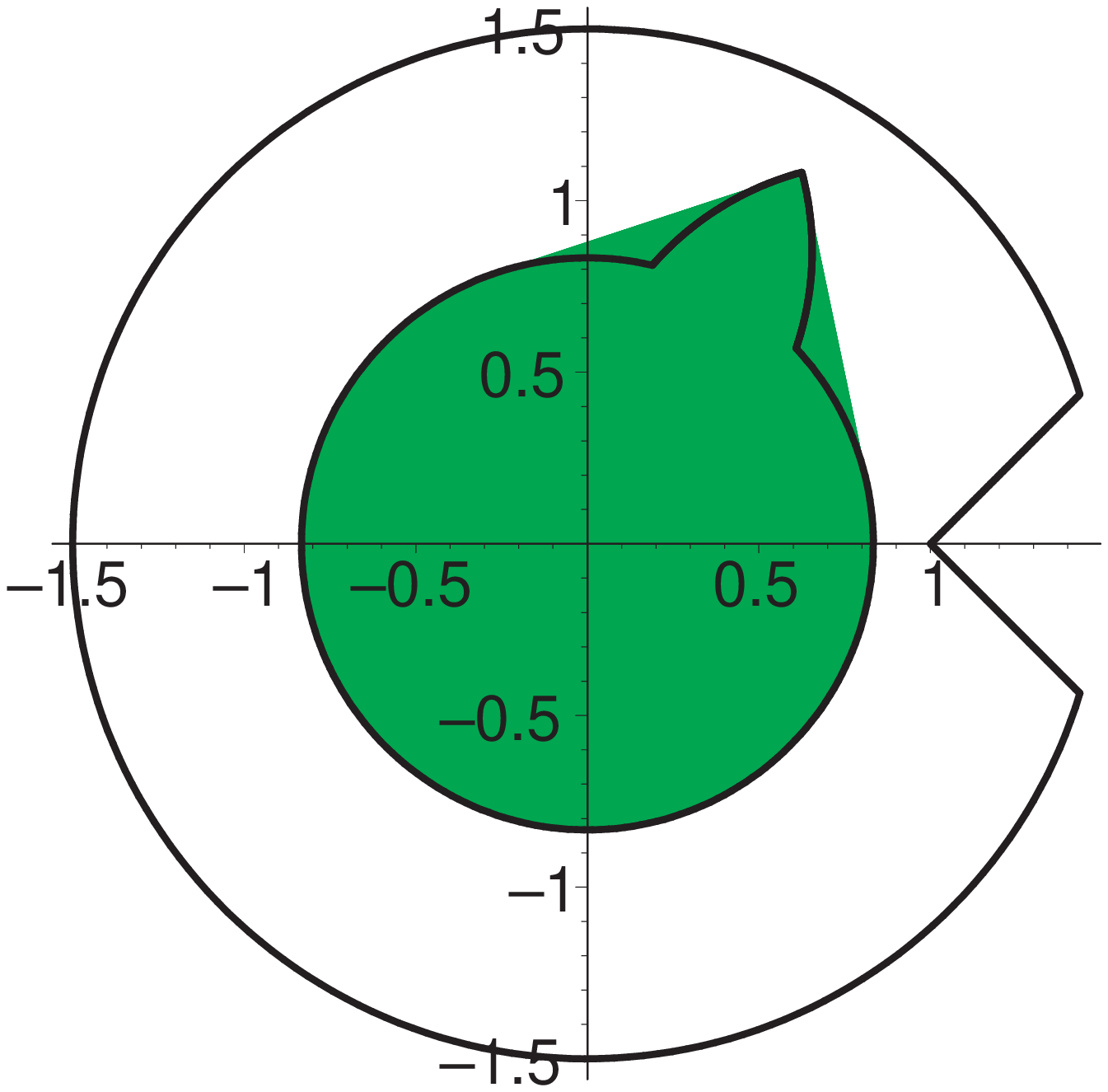}
\end{center}
\caption{\label{hadgeom-fig} The geometry of Hadamard domains:
(left) boundary of a $\Delta$-domain ($1+\eta=1.5$);
(middle) boundary of $\Delta^{-1}$;
(right) an allowable domain in \mbox{$\Delta\cap z\Delta^{-1}$} for
application of Hadamard's formula is the unshaded subset of~$\Delta$
($|z|=1.25$).}%
\end{figure}

\begin{proof} We first observe\footnote{%
    This part of the argument is an adaptation to our needs of a
    famous result first due to Hadamard regarding the continuation of
    Hadamard products; see for instance the description
    in~\cite[Vol.~II, p.~300]{Bieberbach31}
    or~\cite[Sec.~88]{Dienes58}.  Accordingly, we limit ourselves to a
    succinct discussion only meant to set the stage for the precise
    estimates starting at~(\ref{condh2}).} that $f \odot g$
  is continuable to \emph{certain} points $z$ such that $|z| > 1$.
  (Precisely, as shown below, it admits a continuation in a
  $\Delta$-domain.)  Indeed, because of the analytic continuation
  properties of $f$ and $g$, both~$f(w)$ and~$g(z/w)$ are analytic
  functions of~$w$ in the domain $\Delta\cap(z\Delta^{-1})$,
   where $\Delta^{-1}$ denotes $\{w^{-1}: w\in\Delta\}$;
  see Figure~\ref{hadgeom-fig} for a rendering.  In other words, the
  allowed domain of values of~$w$ is $\Delta$ stripped of the internal
  domain $(z\Delta^{-1})^c$, where
    $(\cdot)^c$ represents 
  complementation.  Fix then some~$z_1$ outside the unit disc but
  within~$\Delta$, and choose a simple contour~$\gamma_1$ inside both
  $\Delta$ and $z_1\Delta^{-1}$. Let $I(z)$ be the integral
  of~(\ref{eq:45}) and~(\ref{eq:45bis}) taken along this fixed
  contour~$\gamma_1$. (The feasibility of finding a
  suitable~$\gamma_1$ is suggested by Figure~\ref{hadgeom-fig}, at
  least when~$|z_1|$ remains close enough to~1 and
$z_1$ is to the left of~1; 
a particular contour adapted to the case where~$z_1$ is close to~1 and
possibly to its right will be constructed explicitly  in the proof
below.)  Now, when $z$ moves radially along the segment $(0,z_1)$, the
quantity~$I(z)$ defines an analytic function of~$z$ that does coincide
with~$(f\odot g)(z)$ as soon as $|z|\le 1$  [this results from the
``standard'' formula~(\ref{eq:45bis})].  Thus analytic continuation of
$f\odot g$, from within the unit disc to some~$z_1$ lying outside of
the unit disc
is granted.  The argument shows at the same time that Hadamard's
formula~(\ref{eq:45bis}) remains a valid representation of~$f\odot g$
along such a contour~$\gamma_1$ or any of its deformations legally
granted by analyticity.

\smallskip

%


We next turn to estimating the growth at its singularity of~$h:=f
\odot g$.  It suffices to prove the estimate~(\ref{eq:44}) on $h$
for~$z$ belonging to a restricted domain
$\Delta':=\Delta(\psi_1,\eta_1)$, where we shall take
\begin{equation}\label{condh2}
\eta_1=c_1\eta,\qquad \left(\frac\pi2-\psi_1\right)
=c_1\left(\frac\pi2-\psi_0\right),
\end{equation}
for some small positive constant~$c_1$. Notice also that it suffices
to establish the estimate of~(\ref{eq:44}) for
\begin{equation}\label{condh22}
|z-1|<\eta_1=c_1\eta
\end{equation}
with~$z\in\Delta'$, since~$h$, being analytic in the rest
of~$\Delta'$, is certainly bounded there.

\begin{figure}
  \centering
  \includegraphics[width=7.5truecm]{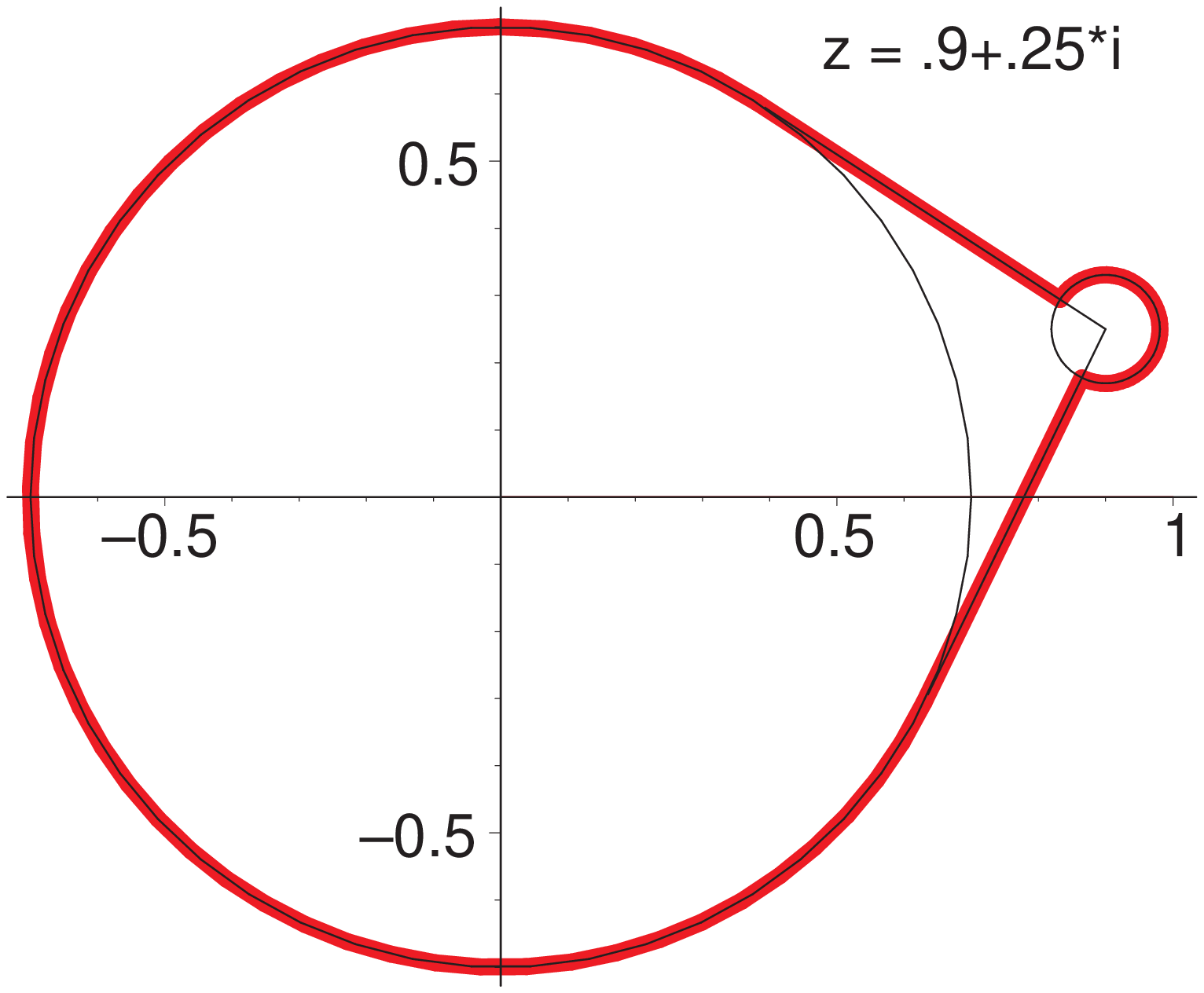}\hspace*{-1.75truecm}
  \includegraphics[width=7.5truecm]{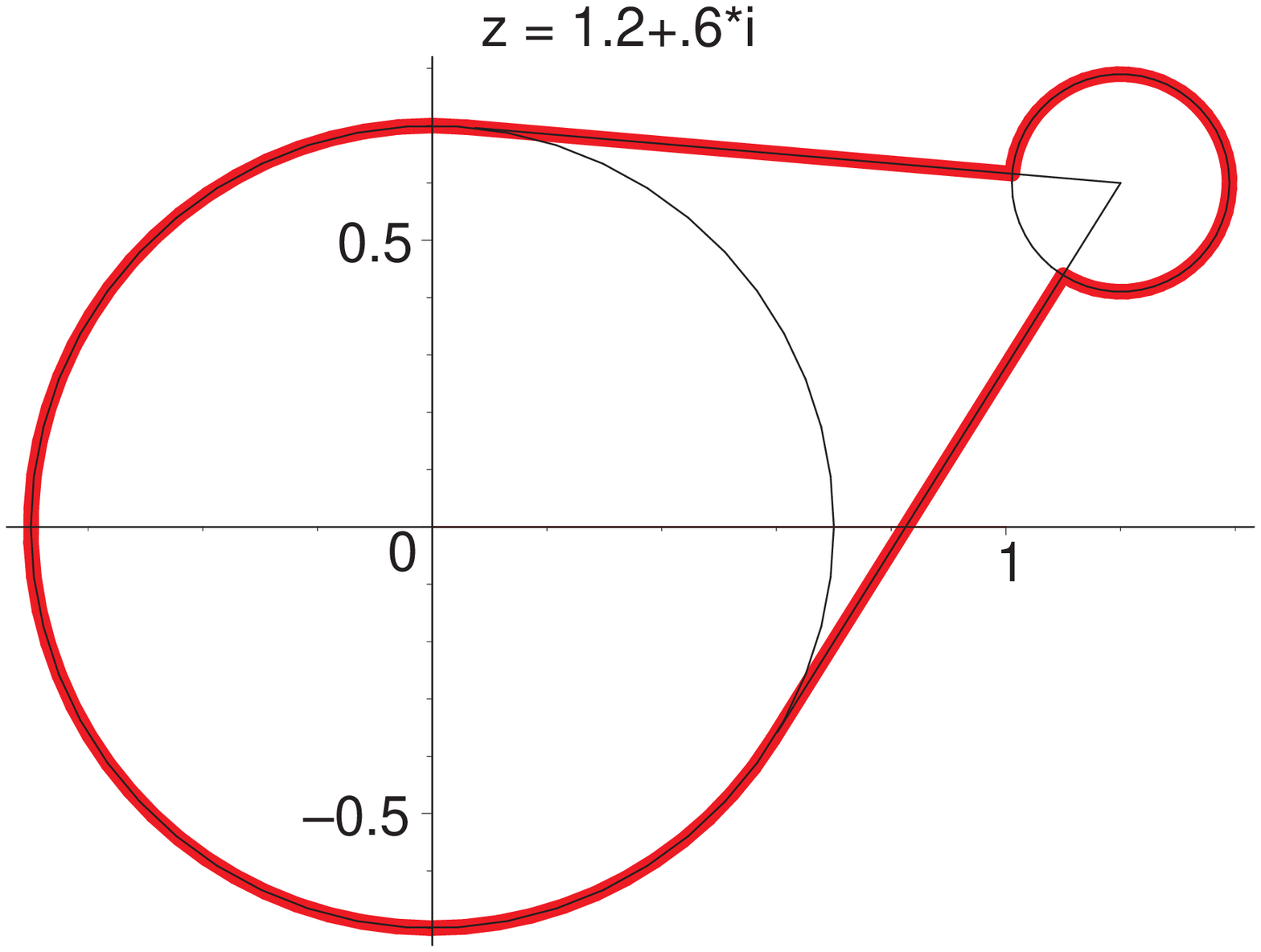}

  \caption{The geometry of the contour $\gamma$.}%
   \label{fig:5-1c}%
\end{figure}

The main geometric objects from which the contour is built are as
follows.  First consider the circle centered at the origin
\begin{equation}\label{condh4}
C_0:=\left\{w: |w|=R\right\},\qquad R:=1-c_2\eta
\end{equation}
for some small constant~$c_2$ (independent of~$z$).
Set~$\delta=|z-1|$,  which is the main parameter
governing the scaling of the contour~$\gamma$. We also consider the
circle 
\begin{equation}\label{condh5}
C_z:=\left\{  w: |w-z| = c_3\delta \right\},
\end{equation}
for some small positive~$c_3$.
Finally, the contour~$\gamma$ includes parts of the two tangents
$T,T'$ to the circle~$C_0$ issuing from~$z$; see
Figure~\ref{fig:5-1c}.  The contour is then precisely specified as
\[
\gamma=\gamma_0 \cup \gamma_T \cup \gamma_z \cup \gamma_{T'},
\]
where $\gamma_T$ is the segment of~$T$ formed of points in between
$C_0$  and~$z$ that are exterior to~$C_0$ and~$C_z$, and
similarly for $\gamma_{T'}$. The component~$\gamma_0$ is the part of
the circle~$C_0$ that lies on the ``southwest'' of~0 and joins
with~$T,T'$; the component~$\gamma_z$ is the part of the circle~$C_z$
that lies on the ``northeast'' of~$z$ and joins with~$T$, $T'$.  The
constants~$c_1,c_2,c_3$ are to be specified  later and
they can be taken as small as needed.

The fundamental constraint to be satisfied is that~$\gamma$ should lie
entirely within~$\Delta\cap(z\Delta^{-1})$ when~$z$ stays
within~$\Delta'$: for~$w\in\gamma$, this ensures simultaneously
$w\in\Delta$ and $z/w\in\Delta$, hence the validity of the Hadamard
integral~(\ref{eq:45bis}).  By a priori  choosing~$c_1$ (which
limits~$z$) and~$c_3$ (which controls the radius of~$C_z$) both small
enough, the condition $\gamma\subseteq\Delta$  is granted
by elementary geometry.  (E.g., the circle $C_z$ will not extend too much
to the right of $\Re(w)=1$ and will therefore be ``compatible'' with
the indentation of~$\Delta$ at~1.)  Next, one should have
$\gamma\cap{(z\Delta^{-1})^c}=\emptyset$.  This requires in
particular choosing the radius~$R$ in~(\ref{condh4}) larger than $z(1
+ \eta)^{-1}$, which 
is at most $(1+c_1\eta)/(1+\eta)$ since $z$ has been restricted to
$|z-1|<c_1\eta$ by~(\ref{condh22}).  This geometric condition
expressed as
\begin{equation}\label{condh6}
\frac{1+c_1\eta}{1+\eta}<1-c_2\eta
\end{equation}
is granted as soon as $c_1,c_2$ are both taken small enough.  [E.g.,
it suffices that both $c_1$, $c_2$ be less than
$\frac12(1+\eta)^{-1}$.]   We
henceforth assume these smallness conditions on~$c_1,c_2,c_3$ to be
satisfied.
Finally, the contour should avoid the apex\footnote{%
  By the ``apex'' of~$(z\Delta^{-1})^c$, we mean the
  complement in~$(z\Delta^{-1})^c$ of the largest circular
  disc centered at the origin which is contained in
  $(z\Delta^{-1})^c$.  } of the domain $(z\Delta^{-1})^c$.
Define the ``viewing angle'' of a point~$P$ exterior to a circle~$C$
as the angle betwen the two tangents to~$C$ issuing from~$P$.  For a
circle of radius~$r$ and a point at distance~$d$ from the center, this
angle is $2\arcsin(r/d)$. In particular the point~$z$ itself views the
circle~$C_0$ of radius~$R$ under the angle $2\arcsin({R}/{|z|})$, and this
viewing angle is bounded from below by
\[
2\arcsin \left(\frac{R}{1+c_1\eta}\right)=2\arcsin\left(
  \frac{1-c_2\eta}{1+c_1\eta}\right),
\]
since the farthest~$z$ can get from the origin is by
assumption~$1+c_1\eta$.  It then suffices to choose~$c_1,c_2$ so that
\begin{equation}\label{view}
2\arcsin\left(
\frac{1-c_2\eta}{1+c_1\eta}\right)>2\psi_0
\end{equation}
(e.g., decide $c_2=c_1$, then decrease~$c_1=c_2$ until the inequality
in~(\ref{view}) is satisfied) in order to ensure that the angle under
which~$z$ views the circle~$C_0$ exceeds~$2\psi_0$.  Since~$C_0$
encloses  the inner disc of~$(z\Delta^{-1})^c$ with
which it is concentric, and since the angle at~$z$ of the apex 
of~$(z\Delta^{-1})^c$ is $2\psi_0$, there results that the
angle at~$z$ between~$\gamma_{T}$ and~$\gamma_{T'}$ encompasses the
apex of~$(z\Delta^{-1})^c$; see Figure~\ref{geomangle-fig}.
In this way, the apex of~$(z\Delta^{-1})^c$ is avoided.

\begin{figure}

\begin{center}
  \includegraphics[width=4truecm]{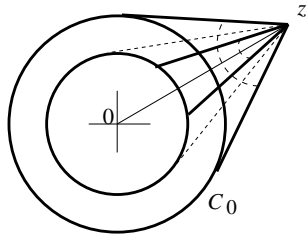}
\end{center}

\caption{\label{geomangle-fig}
  Apex avoidance condition: The angle at~$z$ of contour~$\gamma$ is
  constructed to
  be larger than the angle at~$z$ of the apex of 
  $(z\Delta^{-1})^c$.}
\end{figure}

Last, for~$\lambda$ any of the four contours of which~$\gamma$ is
comprised, let $I(\lambda)$ be the integral of~(\ref{eq:45}) taken
along contour~$\lambda$.  The circular arc $\gamma_{z}$ 
has all its points at a distance  $c_3\delta$ from~$z$,
so that there
\[
|1-w|=\Theta(\delta),\quad |z-w|=\Theta(\delta),\quad
f(w)g\left(\frac{z}{w}\right)=O\left(\delta^{a+b}\right).
\]
Therefore, by trivial bounds,
\begin{equation}\label{hfin1}
I({\gamma_{z}})=O\left(\delta^{a+b+1}\right).
\end{equation}
On the other hand,  along $\gamma_{0}$ the functions
$f(w)$ and $g(z/w)$ stay away from their singularities, so that
\begin{equation}\label{hfin2}
I({\gamma_{0}})=O(1).
\end{equation}
There remains only  to estimate the contribution along
the two connecting segments~$\gamma_{T}$ and $\gamma_{T'}$.  The two
situations are similar (upon interchanging the roles of~$a$ and~$b$).
 It is then easily seen that the contribution along the
ray stemming from $z$ is bounded from above by a multiple of an
integral of the form
\begin{equation}\label{anint}
\int_{c_3\delta}^{+\infty} t^a |t-z_0|^b\, dt
\end{equation}
where~$z_0$ is a complex number at a distance $\Theta(\delta)$ from
the real line.  (The quantity~$t$ parameterizes the tangent line~$T$
or~$T'$.)   The last integral is $O(\delta^{a+b+1})$ as
results from the change of variables~$t=\delta\tau$.  Consequently,
one finds
\begin{equation}\label{hfin3}
I({\gamma_{T}})+I({\gamma_{T'}})=O\left(\delta^{a+b+1}\right).
\end{equation}
Putting together all the estimates of~(\ref{hfin1}), (\ref{hfin2}),
(\ref{hfin3}) yields the desired result. 
\end{proof}

\begin{remark*} 
  \label{rem:intlog}
  The proof technique of Proposition~\ref{lem:3}
tolerates the presence of logarithmic factors,
in which case it suffices to
develop the corresponding estimates for the basic
integral~(\ref{anint}). We find in this way,
when $a+b+1<0$, the estimate 
\[
O\left(|1-z|^a |L^k(z)| \right)\odot O\left(|1-z|^b |L^\ell(z)| \right)
=O\left(|1-z|^{a+b+1} |L^{k+\ell}(z)| \right).
\]
The contour $\gamma$ used in the proof 
is also susceptible to many variations. For instance, one may deform it
slightly to include a ``hook'' near $w=1$,
in which case the modified contour may be used to estimate 
more finely the singular behaviour of Hadamard products.
\end{remark*}
We can then extend the asymptotic range
covered by Proposition~\ref{lem:3} as follows.
\begin{proposition}   \label{thm:4}
   Assume that $f(z)$ and $g(z)$ are $\Delta$-regular and that
for~$z\in\Delta$, 
\[
   f(z)  = O(|1-z|^a)\mbox{\rm \ \ and\ \ }g(z) = O(|1-z|^b).
\] 
   \textit{(i)} If $k < a + b +1 < k+1$ for some integer $-1 \leq k < 
     \infty$, then for~$z\in\Delta'$:
     \begin{equation*}
       \label{eq:46}
     (f \odot g)(z) = \sum_{j=0}^k \frac{(-1)^j}{j!}
     (f\odot{}g)^{(j)}(1)(1-z)^j + O(|1-z|^{a+b+1}).
     \end{equation*}
     \par\noindent
   $(ii)$ If $a+b+1$ is a nonnegative integer then for~$z\in\Delta'$:
     \begin{equation*}
       \label{eq:47}
     (f \odot g)(z) = \sum_{j=0}^{a+b} \frac{(-1)^j}{j!}
     (f\odot{}g)^{(j)}(1)(1-z)^j 
      + O( |1-z|^{a+b+1}|L(z)|).
     \end{equation*}
\end{proposition}
\begin{proof}
Let $\partial=\partial_z$ denote the operator $\tfrac{d}{dz}$ and let $\vartheta$ denote
  the Euler operator $z\partial$. Observe that
     \[
     \vartheta(f\odot{}g) = (\vartheta{}f) \odot{}g =
f\odot (\vartheta g),
     \]
which yields
\[
\vartheta^{k+1}(f\odot{}g)=(\vartheta^{k+1} f)\odot{}g.
\]
The differentiation properties of Theorem~\ref{thm:2A} imply [with $k
:= a + b + 1$ in Case~(ii)] that 
$\vartheta^{k+1}f(z)$ is $O(|1-z|^{a-k-1})$. Thus, Proposition~\ref{lem:3}
applies, to the effect that 
\[
\left( \vartheta^{k+1}(f\odot{}g) \right)(z) = O\left(|1-z|^{a+b-k}\right).
\]
\smallskip
On the other hand, the operator~$\vartheta^{-1}$ is (for $h$ in the
image of $\vartheta$) 
\[
\left( \vartheta^{-1} h \right)(z) :=P_0+\int_0^z h(t)\, \frac{dt}{t},
\]
for some integration constant~$P_0$.
It is then possible to recover~$h=f\odot g$ through successive
integrations,
by making use of Theorem~\ref{thm:2B}.

\smallskip

Case~\textit{(i)}. By definition of $k$,  one has
$-1<a+b-k<0$. Repeated integrations then show that 
\begin{equation}\label{hback1}
(f\odot{}g)(z)=P(z)+O\left(|1-z|^{a+b+1}\right),
\end{equation}
for some polynomial~$P(z)$ of degree~$k$ that encapsulates
the sequence of integration constants. 
Equation~(\ref{hback1})
yields qualitatively the form of the statement. 
The polynomial $P(z)$ is then 
automatically determined as  the first $(k+1)$ terms of
the Taylor expansion of~$f\odot g$ at~$1$, which is precisely what
our assertion  expresses.

\smallskip

Case~\textit{(ii)}. In this case, the first integration step requires
integrating a term~$O(|1-z|^{-1})$, which leads to the 
logarithmic form of the statement. (See also the comments following
Theorem~\ref{thm:2B}.) 
\end{proof}

\subsection{Composition rules}\label{subsec:had-all}
At this stage, we can summarize the state of affairs regarding
Hadamard products by the following general statement.

\begin{theorem}[Hadamard composition of singularities]\label{thm:hadmain}
Let $f(z)$ and $g(z)$ be two functions that are $\Delta$--regular 
with expansions of the type~\eqref{eq:27}: 
\[
f(z)=\sum_{m=0}^M c_m(1-z)^{\alpha_m}+O(|1-z|^A),
\qquad
g(z)=\sum_{n=0}^N d_n(1-z)^{\beta_n}+O(|1-z|^B).
\]
Then, the Hadamard product $(f\odot g)(z)$ is  
also $\Delta$--regular. Its singular expansion is computable 
by bilinearity, using the composition rules of Proposition~\ref{thm:3}
and the remarks thereafter, 
with error terms provided by Propositions~\ref{lem:3} and~\ref{thm:4}:
\[
(f\odot g)(z)=\sum_{m,n} c_m d_n [(1-z)^{\alpha_m}\odot(1-z)^{\beta_n}]
+P(1-z)+
O\left(|1-z|^{C}\right),
\]
where~$C := 1+\min(\alpha_0+B,A + \beta_0)$ 
and~$P$ is a polynomial of degree less than~$C$.
\end{theorem}

The polynomial~$P$ is accessible via the Taylor expansion of
$h-h_{\operatorname{sing}}$, where $h_{\operatorname{sing}}$
represents the sum of all the elements in the asymptotic expansion
of~$h := f \odot g$  at~$z=1$ that are singular.
This theorem then validates the following  algorithm,
which is often helpful in
computations done by hand when composing functions 
under Hadamard products.
\begin{itemize}\def\Asympt{\operatorname{Asympt}}
\item[]
\emph{\bf ``Zigzag'' Algorithm}. [Computes the singular expansion of
$f\odot g$ up to $O(|1-z|^C)$. ]\\[1truemm] 
{\bf 1.} Use singularity analysis to determine 
separately the asymptotic
expansions $\Asympt(f_n)$, $\Asympt(g_n)$ of~$f_n=[z^n]f(z)$
and~$g_n=[z^n]g(z)$ into descending powers of~$n$.
\\
{\bf 2.} Perform the resulting product and reorganize
as~$\Asympt(f_ng_n)$.
\\
{\bf 3.} Choose a basis $\mathcal{B}$ of singular functions,
for instance, the standard basis 
$\mathcal{B}=\left\{(1-z)^\alpha L(z)^k\right\}$,
or the polylogarithm basis,
$\mathcal{B}=\left\{ \Li_{\beta,k}(z)\right\}$.
Construct a function~$H(z)$ expressed
in terms of~$\mathcal{B}$ whose singular behaviour is such that
the asymptotic form of its coefficients, $\Asympt(H_n)$, is compatible
with~$\Asympt(f_ng_n)$ up to the needed error terms.
\\
{\bf 4.} Output the singular expansion of 
$f\odot g$ as the quantity $H(z)+P(z)+O(|1-z|^C)$, where $P$ is a
polynomial in~$(1-z)$ of degree less than~$C$.
\end{itemize}
The reason for the addition of a polynomial in Step~{\bf 4},
is that integral powers of~$(1-z)$ do not leave a trace
in coefficient asymptotics since their contribution is 
asymptotically null. 
(An example of such ``hidden''
analytic terms already appears 
in the composition rule for powers given in Proposition~\ref{thm:3}.)
 The Zigzag Algorithm is  then principally useful for
determining the divergent part of expansions. 
If needed, the coefficients
in the polynomial~$P$ 
can be expressed as 
values of the function $f\odot g$ and
its derivatives at~1 once it has been 
stripped of its nondifferentiable terms.
(This is analogous to the situation prevailing in Proposition~\ref{thm:4}.)

\begin{example} {\em The return of P\'olya's drunkard.}
In the $d$--dimensional lattice~$\mathbb{Z}^d$ of points with integer
coordinates, the drunkard performs a random walk starting from the
origin with steps in~$\{-1,+1\}^d$, each taken with equal
likelihood. The probability that the drunkard is back at the 
origin after $2n$ steps is
\begin{equation}\label{qdrunk}
q_n^{(d)}=\left(\frac{1}{2^{2n}}\binom{2n}{n}\right)^d,
\end{equation}
since the walk is a product~$d$ independent
1--dimensional walks. The probability that $2n$ is the epoch of the
\emph{first} return to the origin is the quantity~$p_{n}^{(d)}$,
which is determined implicitly by
\begin{equation}\label{pdrunk}
\left(1-\sum_{n=1}^\infty p_{n}^{(d)}z^n\right)^{-1}
=\sum_{n=0}^\infty q_{n}^{(d)}z^n,
\end{equation}
as results from the convolution equations
expressing the decomposition of loops into primitive loops.
In terms of the associated ordinary generating functions $P$ and~$Q$, this
relation thus reads as $(1-P(z))^{-1}=Q(z)$. 

The asymptotic analysis of the $q_n$'s is straightforward;
 the one of
the $p_n$'s is more involved and is of interest in connection with
recurrence and transience of the random walk;
see, e.g.,~\cite{DoSn84,Lawler91}.
The Hadamard closure theorem provides a
direct access to this problem. Define
\[
\lambda(z):=\sum_{n\ge0}\frac{1}{2^{2n}}\binom{2n}{n}z^n
\equiv \frac{1}{\sqrt{1-z}}.
\]
Then, Equations~(\ref{qdrunk}) and~(\ref{pdrunk}) imply: 
\[
P(z)=1-\frac{1}{\lambda(z)^{\odot d}}, \qquad
\hbox{where}\quad
\lambda(z)^{\odot d}:=\lambda(z)\odot\cdots \odot\lambda(z)
\ \hbox{($d$ times)}.
\]
The singularities of $P(z)$ are found to be
as follows.

\smallskip
\underline{$d=1$}: No Hadamard product is involved and 
\[
P(z)=1-\sqrt{1-z},\qquad\hbox{implying}\quad
p_n^{(1)}=\frac{1}{n2^{2n-1}}\binom{2n-2}{n-1}\sim\frac{1}{2\sqrt{\pi
n^3}}.
\]
(This agrees with the classical 
combinatorial solution expressed in terms of Catalan numbers.)

\smallskip
\underline{$d=2$}: By the Hadamard closure theorem, the function
$Q(z)=\lambda(z)\odot\lambda(z)$ admits a priori 
a singular expansion at~$z=1$ that is 
composed solely of elements of the form $(1-z)^\alpha$ possibly
multiplied by integral
powers of the logarithmic function~$L(z)$. From a computational
standpoint (cf.\ the Zigzag Algorithm), 
it is then best to start from the coefficients themselves,
\[
q_n^{(2)}\sim \left(\frac{1}{\sqrt{\pi n}}-\frac{1}{8\sqrt{\pi
n^3}}+\cdots\right)^2
\sim
\frac{1}{\pi}\left(\frac{1}{n}-\frac{1}{4n^2}+\cdots\right),
\]
and reconstruct the only singular expansion that is compatible,
namely 
\[
Q(z)=\frac{1}{\pi}L(z)+K+O((1-z)^{1-\epsilon}),
\]
where~$\epsilon > 0$  is an arbitrarily small constant
and $K$ is fully determined as the limit as $z\to1$ of $Q(z)-\pi^{-1}
L(z)$.  Then it can be seen that the function~$P$ is
$\Delta$--continuable.  (Proof: Otherwise,  there would
be complex poles arising from zeros of the function~$Q$ on the unit
disc, and this would entail
in~$p_n^{(2)}$ the presence of terms
oscillating around~0, a fact that contradicts the necessary positivity
of probabilities.)  The singular expansion of~$P(z)$ at~$z=1$ results
immediately from that of~$Q(z)$: 
\[
P(z)\sim 1-\frac{\pi}{L(z)}+\frac{\pi^2 K}{L^2(z)}+\cdots.
\]
so that, by the extension of Theorem~\ref{thm:basic-sa} to arbitrary
powers of logarithms as given in~\cite{FlOd90b,Odlyzko95}, one has
\[\begin{array}{lll}
p_n^{(2)} &=&\ds  \frac{\pi}{n\log^2
n}
-2\pi\frac{\gamma+\pi K}{n\log^3 n}
+O\left(\frac{1}{n\log^4n}\right)
\\
K&=&\ds 1+\sum_{n=1}^\infty
\left(16^{-n} \binom{2n}{n}^2-\frac{1}{\pi n}\right)\\
&\doteq& 0.8825424006106063735858257\,.
\end{array}
\]
(See the study by Louchard \emph{et al.}~\cite[Sec.~4]{LoScToZi94} for
somewhat similar calculations.)

\smallskip
$d=3$: This case is easy since $Q(z)$ remains finite at its
singularity $z=1$ where it admits an expansion in powers of
$(1-z)^{1/2}$, to the effect that
\[
q_n^{(3)}\sim \left(\frac{1}{\sqrt{\pi n}}-\frac{1}{8\sqrt{\pi
n^3}}+\cdots\right)^3
\sim \frac{1}{\pi^{3/2}}
\left(\frac{1}{n^{3/2}}-\frac{3}{8n^{5/2}}+\cdots
\right).
\]
The function~$Q(z)$ is a priori $\Delta$-continuable 
and its singular expansion can be reconstructed
from the form of coefficients: 
\[
Q(z)\mathop{\sim}_{z\to1}
Q(1)-\frac{2}{\pi}\sqrt{1-z}+O(|1-z|),\]
leading to 
\[
P(z)=\left(1-\frac{1}{Q(1)}\right)-\frac{2}{\pi Q^2(1)}\sqrt{1-z}
+O(|1-z|).
\]
By singularity analysis, the last expansion gives
\[
\begin{array}{lll}
p_n^{(3)}&=&\ds
\frac{1}{\pi^{3/2}Q^2(1)}\,\frac{1}{n^{3/2}}+O\left(\frac{1}{n^2}\right)
\\
Q(1) &=& \ds\frac{\pi}{\Gamma\left(\frac{3}{4}\right)^4}
~ \doteq~1.3932039296856768591842463.
\end{array}\]
A complete asymptotic expansion in powers $n^{-3/2},n^{-5/2},\ldots$
can be obtained by the same devices.
In particular this improves the error term above to~$O(n^{-5/2})$.
The explicit form of~$Q(1)$
results from its expression as the generalized hypergeometric 
${}_3F_2[\frac12,\frac12,\frac12;1,1;1]$, 
which evaluates by Clausen's theorem and
Kummer's identity to the square of a 
complete elliptic integral.
(See the papers by Larry Glasser for context, for
instance~\cite{Glasser72}; nowadays, 
{\sf Maple} and {\sf Mathematica} even provide this value
automatically). 

Higher dimensions are treated similarly, with logarithmic terms
surfacing in asymptotic expansions
for all even dimensions.~\hfill~$\Box$  
\end{example}

We observe that, without the developments of the present paper,
the precise asymptotic structure of such sequences is not obvious.
Methods of the last example may be used to provide 
a rigorous setting to certain asymptotic enumeration
results stated by physicists, where back-and-forth  equivalences
between singular expansions of functions 
and asymptotic expansions of coefficients are often used without
much justification. 
See for instance the works of Guttmann and
collaborators~\cite{BoGu97,GuPr93} and Chyzak's 
numeric-symbolic study~\cite{Chyzak97}
relative to special self-avoiding  polygons.

\section{Applications: first moments}\label{sec:some-applications}

Thanks to the extended singularity
analysis toolkit,
we are now in a position to analyze the tree recurrences that were introduced
in Section~\ref{sec:some-special-tree}.
For each of the three models, two types of tolls are
to be considered: 
\[
t_n^{\alpha}:=n^\alpha \text{ (with $\alpha > 0$)},\qquad
t_n^{\log}=\log n,
\]
and we assume in both cases that $t_0=0$.
The corresponding ordinary generating functions are the
polylogarithm~$\Li_\alpha\equiv \Li_{\alpha,0}$ and the specific 
$\Li_{0,1}$,  whose singular expansions have been 
already recalled as Theorem~\ref{thm:li-sing}.
In each case, a linear transform~$\cal L$ relates 
the generating function of costs, $f(z)$, to a generating function of
tolls, either $t(z)$ (normalized) or $\tau(z)$ (``raw''). 
Theorems on composition of singularities make it possible to 
follow step
by step the elementary operations of which~$\cal L$ is composed
and determine the effect of the $\cal L$ transform on singularities in
a systematic manner. Given that computations are ``automatic'', we
will mostly focus our discussion on main terms and on the global shape
of singular expansions, leaving some of the details as exercises to
the  reader---or better, to a computer algebra engine.

The net outcome in each of the three tree models under consideration
is the following: for large tolls, the cost is driven by the toll itself;
for small tolls, the cost is of linear growth and, in a sense, 
``freely'' caused by
the recursion itself, that is, driven by the cumulation of costs
due to small subtrees; in between, 
there is a threshold value of the toll where a ``resonance''
takes place between the toll and the recursion, leading to the emergence of a
logarithmic factor. Such facts parallel what is familiar in the context
of inhomogeneous linear differential equations, where either the free
regime or the forced regime dominates,  with logarithmic terms
being created precisely  by resonances.


\subsection{The binary search tree recurrence}\label{sec:binary-search-tree}

For the binary search tree model, 
there is an integral transform $\mathcal L$ that relates the 
ordinary generating function of tolls, $t(z)$, and the ordinary generating
function of the induced costs, $f(z)$: it is given by~(\ref{bst3})
according to which $f(z)=\mathcal L[t(z)]$, where (with~$t_0=f_0=0$)
\begin{equation}\label{Lbst}
\mathcal L[t(z)]=\frac{1}{(1-z)^2}\int_0^z t'(w)(1-w)^2\, dw.
\end{equation}
Consequently, the computation is entirely \emph{mechanical}\footnote{%
        In the {\sc Maple} system for symbolic computations
about two dozen instructions suffice to implement calculations,
once use is made of Bruno Salvy's package {\sf equivalent}
dedicated to the asymptotic analysis of coefficients of generating
functions~\cite{Salvy91}.
It suffices to program the polylogarithm expansions
(Theorem~\ref{thm:li-sing}),
use the system capabilities for series expansions, differentiation, and
integration (Theorems~\ref{thm:2A} and~\ref{thm:2B}), and conclude by
an appeal to Salvy's program that implements the basic transfers of
Theorem~\ref{thm:basic-sa}. 
} and
it only needs the theorems relating to integration, differentiation,
and polylogarithms (Theorems~\ref{thm:2A} and~\ref{thm:2B}) in 
conjunction with basic 
singularity analysis (Theorem~\ref{thm:basic-sa}). 
Our derivation below constitutes an alternative to
parallel results by Neininger~\cite{Neininger02},
Chern~\emph{et al.}~\cite{ChHw01,ChHwTs02}, and  Fill and
Kapur~\cite{FK-transfer}, who employ elementary but
perhaps less transparent
methods (typically, the approximation of discrete sums by
integrals).

\begin{theorem} Under the binary search tree model,
the expected values of the costs 
induced by tolls of type $t_n^\alpha$ 
$(\alpha > 0$) and $t_n^{\log}$ admit full asymptotic expansions
in descending powers of~$n$ and integral powers of~$\log n$.
The main terms are summarized by the following table: 
\[
\begin{array}{cl|cl}
\hline
\hbox{\em Toll $(t_n)$} & & \hbox{\em Cost $(f_n)$} & \\ 
\hline\hline
n^\alpha & (2<\alpha) & \ds \frac{\alpha+1}{\alpha-1}n^\alpha &
{}+O(n^{\alpha-1}) \\[.15in]
n^2 &  & \ds 3n^2 &{}-6n\log n+(10-6\gamma )n+O(\log n) \\[.1in]
n^\alpha & (1<\alpha<2) & \ds \frac{\alpha+1}{\alpha-1}n^\alpha &
{}+K_\alpha n+O(n^{\alpha-1}) \\
n & & 2n\log n &\ds {}+2(\gamma-1)n+2\log n+2\gamma+1+O\left(\frac1n
\right) \\
n^\alpha & (0<\alpha<1) & K_\alpha n &
\ds {}+\frac{\alpha+1}{\alpha-1}n^\alpha+K_\alpha+o(1)\\ 
\log n & & K'_0 n & \ds {}-\log
n+(K'_0-2)-\frac{1}{2n}+\frac{1}{9n^2} +O\left(\frac{1}{n^3}\right).\\ 
\hline
\end{array}
\]
\end{theorem}
\begin{proof} For the case $\alpha$ a nonnegative integer, 
  the integration can be carried out in finite terms since the
  generating function of tolls is rational. For instance, the case
  $\alpha=1$ corresponds to the well-known analysis of 
  Quicksort and binary search tree 
  algorithms~\cite{Knuth98a,Mahmoud92,SeFl96,ViFl90}.
  
  \def\leadsto{\ \longrightarrow\ } For $t_n^\alpha$, it suffices to
  examine the effect of the $\mathcal L$ transform on singular
  elements of the form $c(1-z)^\beta$; e.g., for the main term
  corresponding to $t_n=n^\alpha$, we should take $\beta=-\alpha-1$.
  The $\mathcal L$ transformation reads as a succession of operations,
  ``differentiate, multiply by~$(1-z)^2$, integrate, multiply
  by~$(1-z)^{-2}$''---all are covered by our previous theorems. The
  chain on any particular singular element starts as
\[
c(1-z)^\beta \mathop{\leadsto}^{\partial} c\beta(1-z)^{\beta-1}
\mathop{\leadsto}^{\times(1-z)^2} c\beta(1-z)^{\beta+1}.
\]
At this stage, integration  intervenes. Assume that $\beta+1\not=-1$.
(Otherwise, a logarithm appears.)  According to Theorem~\ref{thm:2B},
and ignoring integration constants for the moment, integration gives
\[
c\beta(1-z)^{\beta+1}\mathop{\leadsto}^{\int}
-c\frac{\beta}{\beta+2}(1-z)^{\beta+2}
\mathop{\leadsto}^{\times(1-z)^{-2}}
-c\frac{\beta}{\beta+2}(1-z)^{\beta} .
\]
Then this singular element corresponds to a contribution
\[
-c\frac{\beta}{\beta+2}\binom{n-\beta-1}{-\beta-1},
\]
which is of order $O(n^{-\beta-1})$.  (The treatment of logarithmic
terms is entirely similar.)

The derivation above
has left aside the determination of the integration constants. These
are given by the second case of Theorem~\ref{thm:2B}, which provides
in particular access to the constants~$K_\alpha$ and $K'_0$. The
constant term in the asymptotic expansion of the integral is of the
form 
\[
{\bf K}[t]:=\int_0^1
\left[t'(w)(1-w)^2-\left(t'(w)(1-w)^2\right)_-\right]\,dw,
\]
where $f_-$ represents the sum of the singular terms in~$f$ having
exponent $<-1$, as in the proof of Theorem~\ref{thm:2B}. In the
singular expansion of~$f(z)$, this integration constant gets further
multiplied by~$(1-z)^{-2}$;  the resulting linear term in the
asymptotic expansion of $f_n$ is then plainly
\[
{\bf K}[t]\cdot (n+1).\] In particular, if the growth of~$t_n$ is
smaller than $n$, then, the divergence part is absent and ${\bf K}[t]$
reduces to
\[
{\bf K}[t]=\int_0^1 t'(w)(1-w)^2\,dw=2\sum_{n=1}^\infty
\frac{t_n}{(n+1)(n+2)},
\]
as follows from expanding the integrand around~0 and integrating the
resulting series.  This yields the following values for $\alpha<1$:
\begin{equation}\label{eq:theKs}
K_\alpha=2\sum_{n=1}^\infty \frac{n^\alpha}{(n+1)(n+2)},
\qquad 
K'_0=2\sum_{n=1}^\infty \frac{\log n}{(n+1)(n+2)},
\end{equation}
while for~$1<\alpha<2$, 
\[
K_\alpha=2\sum_{n=1}^\infty
\frac{n^\alpha-\Gamma(\alpha+1)\binom{n+\alpha}{\alpha}}{(n+1)(n+2)}.
\]
The theorem is finally established.
\end{proof}
\begin{remark*}
The slowly convergent series expressions 
of $K_\alpha,K'_0$ can be rephrased as definite
integrals, thanks to Mellin transform techniques. The starting point
is the easy formal identity, 
\begin{equation}\label{mellin0}
\sum_{n\ge1} c_n n^{-s} = \frac{1}{\Gamma(s)}
\int_0^\infty \left(\sum_{n\ge1} c_n e^{-nx}\right) x^{s-1}\, dx.
\end{equation}
The constant~$K_\alpha$ with~$\alpha<1$ corresponds to~$s=1-\alpha$
 and 
 $c_n=n/[(n+1)(n+2)]$, for which 
the integrand admits of closed form since 
\[
\sum_{n=1}^\infty
\frac{n z^n}{(n+1)(n+2)}= \frac{1}{z^2}\left[(2-z)L(z)
-2z \right].
\]
From there, the constant $K'_0$ is  attained as 
$\left.\frac{d}{d\alpha}K_\alpha\right|_{\alpha=0}$. 
A final change of variables $x=-\log t$ then yields
an integral representation for
``Fill's first logarithmic constant'' ($\gamma$ is Euler's constant):
\begin{eqnarray}\label{fillcbst}
K'_0&=&-\gamma-2\int_0^1\left[(t-2)\log(1-t)-2t\right]
\left(\log\log\frac1t\right) \, \frac{dt}{t^3}  \\
&=&
1.20356491674961033428628333814873131775552838577096
. \nonumber
\end{eqnarray} 
The last estimate to 50D improves on the earlier 3D  evaluation
of~Fill \cite{Fill96}. The cost induced by~$t^{\log}$ is of 
particular interest as
it is precisely the entropy of the distribution of binary search
trees; see the account and first estimates in the book by Cover and
Thomas~\cite[p.~72--74]{CoTh91}, as well as pointers to
self-organizing search in Fill's article~\cite{Fill96}. 
In his doctoral dissertation~\cite[Section~5.1]{kapur03:_addit},
Kapur has extended the methods and estimates
to~$m$-ary search trees. 
\end{remark*}

\subsection{The uniform binary tree recurrence}\label{sec:binary-tree-recurr}

This section examines the uniform binary tree model 
that surfaces recurrently in combinatorics. 
Here, we put on a firm basis a classification of the expected costs
corresponding the tolls $t_n^\alpha$ and $t_n^{\log}$
which was outlined (with several typographical errors) in 
an article by Flajolet and Steyaert~\cite{FlSt87}.
The particular case of the toll~$t_n=n$ has, like for binary search
trees, a dignified history as it corresponds to path length
in Catalan trees and to area under Dyck paths, whose first distributional
analyses go back to Louchard and Tak\'acs~\cite{Louchard84,Takacs91}.

Our starting point is~(\ref{cat99})  according to which the
generating function of costs~$f(z)=\sum C_nf_nz^n$ normalized by the 
Catalan numbers~$C_n$ and the ordinary generating function of costs
$\tau(z)=\sum t_n z^n$ are related by $f(z)=
\overline{\cal L}[\tau (z)]$, where
\begin{equation}\label{Lcat}
\overline{\cal L}[\tau(z)]=\frac{1}{\sqrt{1-4z}}\left(\tau(z)\odot C(z)\right),
\end{equation}
with
\[
C(z)=\sum_{n\ge0} C_n z^n =\frac{1}{2z}\left(1-\sqrt{1-4z}\right),
\qquad
C_n=\frac{1}{n+1}\binom{2n}{n}.
\]
We state:
\def\K{\overline{K}}
\begin{theorem} Under the uniform binary tree model,
the expected values of the costs 
induced by tolls of type $t_n^\alpha$
$(\alpha > 0$)  and $t_n^{\log}$ admit full asymptotic expansions
in descending powers of~$n$ and integral powers of~$\log n$.
The main terms are summarized by the following table: 
\[
\begin{array}{cl|cl}
\hline
\hbox{\em Toll $(t_n)$} & & \hbox{\em Cost $(f_n)$} & \\ 
\hline\hline
n^\alpha & (\frac32 < \alpha) & \ds
\frac{\Gamma(\alpha-\frac12)}{\Gamma(\alpha)}n^{\alpha+\frac12}& 
{}+O(n^{\alpha-\frac12}) \\[.15in]
n^{3/2} & &
\ds \frac{1}{\Gamma(3/2)}n^2 &{}+O(n\log n)
\\[.15in]
n^\alpha & (\frac12<\alpha<\frac32) & \ds
\frac{\Gamma(\alpha-\frac12)}{\Gamma(\alpha)}n^{\alpha+ \frac12}&
\ds {}+O(n) \\[.15in]
n^{1/2} & & \ds \frac{1}{\sqrt{\pi}} n\log n &\ds 
{}+O(n) \\[.15in]
n^\alpha & (0<\alpha<\frac12) & \K_\alpha n &
\ds {}+O(1)
\\[.15in]
\log n & & \K'_0 n & \ds {}-2\sqrt{\pi}n^{1/2}+O(1).
\\
\hline
\end{array}
\]
\end{theorem}
\begin{proof}
For the tolls $t_n^\alpha$, all that is required is to determine the
singular expansion of
\[
\tau(z)\odot C(\frac{z}{4}) = \sum_{n=1}^\infty \frac{n^\alpha}{(n+1)}
\binom{2n}{n}\left(\frac{z}{4}\right)^n.
\]
(For convenience, the singularity has been scaled to~1.)
We use the Zigzag Algorithm 
presented in Subsection~\ref{subsec:had-all}.
The known asymptotic expansion of the Catalan numbers is
\[
4^{-n}C_n\sim\frac{1}{\sqrt{\pi}} {n^{-3/2}}
\left(1-\frac{9}{8n}+\frac{145}{128n^2}-\cdots\right).
\]
Multiply this by $n^\alpha$ to get the expansion of $n^\alpha C_n4^{-n}$.
The terms now involve the scale 
$\{n^{\alpha-\frac32}, n^{\alpha-\frac52},\ldots\}$.
Assume that $\alpha$ is not a half-integer [i.e., \mbox{$\alpha \not\in
(\frac12\Z) \setminus \Z$}]; see below
for the  contrary case. Then
the basis of functions $\mathcal{B}=\{(1-z)^{-\alpha+k+ \frac12}\}$,
 where~$k$ ranges over the integers, has the property that 
the coefficients of its generic element
are~$O(n^{\alpha-k-\frac32})$;  in particular, 
\[
[z^n](1-z)^{-\alpha+
  \frac12}\sim\frac{n^{\alpha-\frac32}}{\Gamma(\alpha-\frac12)} 
\left(1+\frac{(2\alpha-1)(2\alpha-3)}{8n}+\cdots\right).
\]
We can thus find a singular function $H(z)$ whose coefficients 
match asymptotically those of $\tau(z) \odot C(z/4)$, which is of the
form 
\[
H(z)=\frac{\Gamma(\alpha- \frac12)}{\sqrt{\pi}}
(1-z)^{-\alpha+ \frac12}\left(1+c_1(1-z)+c_2(1-z)^2+\cdots\right),
\]
for some effectively computable sequence~$(c_j)$.
The singular expansion of $\tau(z)\odot C(z/4)$  is
then the  sum of the expansion of~$H$ above and of a
power series in~$(1-z)$, call  it $P(z)$, that can be determined 
according to the principles of Section~\ref{sec:hadam-prod-transf}.


The singular expansion of $f(z/4)$ is that of 
$H(z)+P(z)$ divided by~$\sqrt{1-z}$, so that, by transfer, we get
\[
[z^n]f\left(\frac{z}{4}\right)\sim
\frac{\Gamma(\alpha-\frac12)}{\sqrt{\pi}\Gamma(\alpha)}
n^{\alpha-1}\left(1+\frac{c'_1}{n}+\frac{c'_2}{n^2}+\cdots\right)
+[z^n]\frac{P(z)}{\sqrt{1-z}},
\]
for some sequence~$(c'_j)$, where the ``hidden'' analytic part
$P(z)$ arises from the ``hidden'' analytic 
component in $\tau(z)\odot C(z/4)$. 
 After dividing by~$C_n4^{-n}$, one finds finally: 
\begin{equation}\label{zoba}
f_n \sim \frac{\Gamma(\alpha-\frac12)}{\Gamma(\alpha)}
n^{\alpha+\frac12}\left(1+\frac{c''_1}{n}+\cdots\right)+R_n
,
\end{equation}
where the ``hidden'' remainder term $R_n$ is of the form
\[
R_n\sim d_{-1}n+d_0+\frac{d_1}{n}+\cdots\,.
\]
This last estimate provides all the entries in the table above,
whenever~$\alpha$ is not a half-integer,  as it
suffices to merge the two expansions of~(\ref{zoba}).  In addition, when
$0<\alpha<\frac12$, the series defining $t(z/4)$ converges at the
singularity~1. Thus, the dominant asymptotic term of $f(z/4)$ is
$t(1/4)/\sqrt{1-z}$, that is,
\[
f\left(\frac z4\right)\sim \frac{\K_\alpha}{\sqrt{1-z}},\qquad
\K_\alpha:=\sum_{n=1}^\infty \frac{n^\alpha}{n+1}
\frac{1}{4^n}\binom{2n}{n}.
\]

When $\alpha$ is a half-integer, logarithmic terms appear due to the
presence of inverse integral powers of~$n$ in the coefficients
of $t(z/4)$, but the derivation is otherwise similar.
For instance at $\alpha=\frac12$, one has 
\[
4^{-n}\sqrt{n}C_n\sim \frac{1}{\sqrt{\pi}}\,\frac{1}{n}
+O\left(\frac{1}{n^2}\right),
\]
which shows that 
\[
t(z/4)=H(z)+P_0+O((1-z)^{1-\epsilon}),
\qquad
H(z)=\frac{1}{\sqrt{\pi}}L(z),\]
resulting in the stated estimate.

Finally, when $t_n = \log{n}$, we have $\tau(z) = \Li_{0,1}(z) =
O(|1-z|^{-1-\epsilon})$ for any $\epsilon > 0$. 
Thus, by Proposition~\ref{thm:4}(i), 
\begin{equation*}
  (\tau \odot C)(z/4) = \K'_0+  O\Bigl(|1-z|^{\tfrac12 - \epsilon}\Bigr),
\qquad
  \K'_0 := \sum_{n=1}^\infty (\log n) \frac{C_n}{4^n}.
\end{equation*}
Singularity analysis and the estimate for $C_n$  
yield $f_n = \K'_0 n + O\Bigl(n^{\tfrac12 + \epsilon}\Bigr)$.
Carrying higher-order terms, we get the
mean of the shape functional, 
\begin{equation}
  \label{eq:24}
  \mu_n = \K'_0 n - 2\sqrt\pi n^{1/2} + O(1),
\end{equation}
which agrees with the estimate in Theorem~3.1 of~\cite{Fill96}
and improves the remainder estimate.
\end{proof}

The Mellin technique of~(\ref{mellin0}) is once more applicable
to the determination of ``Fill's
second logarithmic constant''~$\K'_0$. It provides
the value: 
\[\begin{array}{ccc}
\K'_0&:=&\ds \sum_{k\ge1} \frac{\log k}{(k+1){4^k}}\binom{2k}{k}.
\\
&=& \ds -\gamma-\int_0^1 \frac{1}{\sqrt{1-t}(1+\sqrt{1-t})^2}\,
\left(\log\log\frac1t\right)\, dt,\\
\\
&=& 2.0254384677765738877135187391417652470652930617658.
\end{array}
\]
The subject of costs on binary trees is considered in greater depth
in~\cite{fill03:_limit_catal} by
applying the techniques developed in this paper.  There, some higher-order
estimates,
asymptotics for moments of each order, and limiting distributions are
derived when the toll sequence is either $n^\alpha$ or
$\log{n}$.

Our methods can also be used to treat more generally the case of
all simple families of trees in the sense of Meir and
Moon~\cite{MeMo78}, of which Catalan trees are a
special case.  This generalization is the subject of ongoing
work.

\subsection{The union--find tree recurrence}  
\label{sec:cayl-tree-recurr}

\def\K{\widehat{K}}

In this subsection, we revisit the Knuth--Pittel--Sch{\"o}nhage
recurrence corresponding to the destruction of free labelled trees and
dually to the management of equivalence
relations~\cite{KnPi89,KnSc78}.  The main result of this section is
essentially a rephrasing of the main results of Knuth and Pittel
in~\cite{KnPi89}, to which we add the possibility of determining
complete asymptotic expansions.  Like before, the starting point is
the integral transform~(\ref{cay3}) (adjusted for the fact that $t_1^\alpha =
f_1^\alpha = 1 \ne 0$),  which relates the ordinary generating function of
tolls~$\tau(z)$ to the normalized generating function of costs~$f(z)$
via $f(z)=\cal L[\tau(z)]$, where 
\begin{equation}\label{Lcay}
\cal L [\tau(z)]   = t_1 z T'(z) + 
\frac12\frac{T(z)}{1-T(z)}\int_0^z \partial_w
\left(\tau(w)\odot
T^2(w)\right) \, \frac{dw}{T(w)}.
\end{equation}
There~$T(z)$ is the Cayley tree function whose  
singular expansion at the (unique) dominant 
singularity $z=e^{-1}$ is well known:  one has the shape 
\begin{equation}\label{eq:10} 
  T(z) \sim 1 - \sqrt{2}(1-ez)^{1/2} + c_1(1-ez) + \cdots
\end{equation}
as $z \to e^{-1}$ in any sector of angle $<2\pi$; see
also~\cite[Eq.~(3.16)]{KnPi89}.  (The paper by Corless \emph{et
  al.}~\cite{CoGoHaJeKn96} is a definitive reference regarding the
tree function.)  As noted earlier, the case of union--find tree
recurrences 
 combines all the composition results developed in this
paper.
\begin{theorem} Under the union--find tree recurrence model, 
the expected values of the costs 
induced by tolls of type $t_n^\alpha$ 
$(\alpha > 0$) and $t_n^{\log}$ admit full asymptotic expansions
in descending powers of~$n$ and integral powers of~$\log n$.
The main terms are summarized by the following table:  
\[
\begin{array}{cl|cl}
\hline
\hbox{\em Toll $(t_n)$} & & \hbox{\em Cost $(f_n)$} & \\ 
\hline\hline
n^\alpha & (\frac32 < \alpha) & \ds
\frac{\Gamma(\alpha-\frac12)}{\sqrt2\Gamma(\alpha)}n^{\alpha+\frac12}&
{}+O(n^{\alpha-\frac12}) \\[.15in]
n^{3/2} & &
\ds \frac{1}{\sqrt2 \Gamma(3/2)}n^2 &{}+O(n\log n)
\\[.15in]
n^\alpha & (\frac12<\alpha<\frac32) & \ds
\frac{\Gamma(\alpha-\frac12)}{\sqrt 2\Gamma(\alpha)}n^{\alpha+\frac12}&
\ds {}+O(n) \\[.15in]
n^{1/2} & & \ds \frac{1}{\sqrt{2\pi}} n\log n &\ds 
{}+O(n) \\[.15in]
n^\alpha & (0<\alpha<\frac12) & (1 + \frac12\K_\alpha) n &
\ds {}+O(n^{\alpha + \frac12})
\\[.15in]
\log n & & \frac12\K'_0 n & \ds {}+O(\sqrt{n}).
\\
\hline
\end{array}
\]
\end{theorem}

\begin{proof}
\def\A{{\mathcal{A}}}\def\B{{\mathcal{N}}} 
We shall content ourselves with indicating the way
full asymptotic expansions can be determined within the generating
function framework. (Detailed computations are left as an exercise for
the reader.) 
In what follows, we set $Z=(1-z)$ and let~$\A$ 
denote an unspecified entire series in powers of~$Z$, 
not necessarily the same at each
occurrence. 
For instance, one may summarize 
diversely the expansion~(\ref{eq:10}) of~$T(z/e)$ 
as 
\[
T(z/e)\sim 1-\sqrt{2}Z^{1/2}+Z\A+Z^{3/2}A\sim \A+\A Z^{1/2},
\]
and so on. We shall also let $\B$ denote 
generically a series in descending powers
of~$1/n$.

We consider first the case of the toll~$t_n^\alpha$
and assume for simplicity that $\alpha$ is not a half-integer:
$\alpha\not\in(\frac12\Z) \setminus \Z$. 
The polylogarithm expansions grant us 
{a~priori} that the generating
function~$\tau(z)$ lies in the class of functions amenable to
singularity analysis, with 
\[
\tau(z)\sim Z^{-\alpha-1}\A+ \A.
\]
Therefore, the Hadamard product  $(\tau(z)\odot
T^2(z/e))$ is also amenable.   The coefficients of the
latter function are of the form 
$n^{\alpha-\frac32}\B$,  as follows from the fact that
$[z^n]\tau(z)=n^\alpha$ and  $[z^n]T^2(z/e)\sim n^{-3/2}\B$
(by the singular expansion of~$T^2$).  Thus, converting back this
information to the function, we find 
\[
\tau(z)\odot T^2(z/e) \sim Z^{-\alpha+\frac12}\A+\A,
\qquad
\partial_z(\tau(z)\odot T^2(z/e))\sim Z^{-\alpha-\frac12}\A+\A.
\]
What we have done here is to apply  the Zigzag Algorithm of
Section~\ref{sec:hadam-prod-transf} and the differentiation theorem.
Then multiplication by~$1/T(z/e) \sim \A+Z^{1/2}\A$  shows
that 
\[
\frac{1}{T(z/e)}\partial_z[\tau(z)\odot T^2(z/e)]
\sim Z^{-\alpha-\frac12}\A+Z^{-\alpha}\A+\A+Z^{\frac12}\A.
\]
Integration of this last expansion corresponds to
increasing all exponents by~1. Finally one should multiply by
$T(z/e)(1-T(z/e))^{-1}$  which is of type
\mbox{$Z^{-1/2}\A+\A$}.  This completes our handling of the second term
on the right in~(\ref{Lcay}).  Also, 
\[
 \frac{z}e T'(z/e) \sim Z^{-\frac12}\A + \A.
\]
The end result is then 
\[
f(z/e)\sim Z^{-\alpha}\A+Z^{-\alpha+\frac12}\A+Z^{-\frac12}\A+\A.
\]
The dominant term is $Z^{-\alpha}$ when~$\alpha>\frac12$ whereas it is
$Z^{-1/2}$ when $\alpha<\frac12$.

At the same time, it is a simple task to trace the 
coefficients of main terms. For $\alpha>\frac12$, the
main term of~$f(z/e)$ is $Z^{-\alpha}$,  
and one finds successively 
\begin{equation*}\label{eq:cayley18}
\begin{array}{lll}\ds   \tau(z) \odot T^2(z/e)&\sim& 
\ds \sqrt\frac2\pi \Gamma(\alpha -
  \tfrac12) (1-z)^{-\alpha + \frac12},\\[.15in]
  f(z/e)&\sim& \ds\frac{\Gamma(\alpha - \frac12)}{2\sqrt\pi}
(1-z)^{-\alpha},
\end{array}
\end{equation*}
where the last equation implies, via singularity analysis,
an estimate of expected costs: 
\begin{equation*}
  f_n^\alpha \sim  \frac{\Gamma(\alpha - \frac12)}{\sqrt2\Gamma(\alpha)}
  n^{\alpha + \frac12}. 
\end{equation*}

For $\alpha<\frac12$, the main term is $Z^{-1/2}$
and its coefficient is seen to arise from both terms on 
the right in~(\ref{Lcay}):\ we have  
\begin{equation*}
  f^\alpha(z/e) \sim  \frac{1}{\sqrt{2}}\left(1 +
  \frac{\K_\alpha}{2} \right) (1-z)^{-1/2} ,
\end{equation*}
where $\K_\alpha=\K[n^\alpha]$ and the functional~$\K$ is 
\begin{equation}
  \label{eq:16}
  \K[t] := \int_0^{1/e} 
  \partial_w ( \tau(w) \odot T^2(w) )\,  \frac{dw}{T(w)}. 
\end{equation}
Error terms can be similarly traced: 
in the case of~$f^\alpha(z/e)$, it is  
of type~$Z^{-\alpha}$ if~$0<\alpha<\frac12$, of type~$Z^{-1/2}$
if $\frac12<\alpha<1$, and so on.
The end results are summarized in the statement of the theorem.

For half-integer $\alpha$, a logarithmic term appears.
For instance, in the  case $\alpha=\frac12$, this fact is
associated to the shape of the coefficients 
\[
[z^n](\tau(z)\odot T^2(z/e))\sim
\frac{-2\sqrt{2}}{\Gamma(-1/2)}\,\frac{1}{n}
+\frac{1}{n^2}\B,
\]
resulting in a  singular expansion with a logarithmic term:
\[
\tau(z)\odot T^2(z/e)=-\sqrt{\frac{2}{\pi}}\log Z+c+Z\A+(Z\log Z)\A,
\]
for some~$c$.

For the logarithmic toll [note that now $t_1=0$, so that
the first term in~(\ref{Lcay}) does not contribute] we
have~$\tau(z)=\Li_{0,1}(z)$, the integral in~(\ref{Lcay}) is
convergent and, in the same way as for the case $\alpha < 1/2$, we get 
\begin{equation*}
  f(z) = \frac{\K'_0}{2\sqrt2}
  (1-ez)^{-1/2} + O(|1-ez|^{-\epsilon}), \end{equation*}
which implies
\begin{equation*}
  f_n = \frac12\K'_0 n + O(n^{\frac12 + \epsilon}),
\end{equation*}
with   $\K'_0=\K[\log n]$
and $\K$ defined at~\eqref{eq:16}.
\end{proof}

It is of interest to compare our approach
to that of Knuth and Pittel~\cite{KnPi89}.
These authors use what is fundamentally a ``repertoire approach'',
based on the transforms of two types of tolls, the Dirac
tolls~$\delta_{mn}$ and another family related to 
``tree polynomials''. Their methods 
do not clearly appear to be extendible to 
the extraction of sublinear terms in asymptotic expansions.
At the same time, their developments require appreciably more
involved and perhaps less transparent computations.


%

\section{Perspectives}\label{sec:appl-high-moments}

In this concluding section, we discuss at a fairly informal 
and abstract level
applications of the extended
singularity analysis toolkit 
developed in the present paper in two further directions:
the determination of higher-order moments for our 
basic models, and the treatment 
of tree recurrences which are more complex than the ones present in our
lead examples.
(Some of our examples below may accordingly involve
nonbinary tree models.)

\subsection{Higher moments and limit distributions}
Let us return to the general framework of Section~\ref{DandC-sec}. 
There, the random cost~$X_n$ is related to costs~$X_{K_n}$
and~$\widetilde{X}_{n-a-K_n}$
by the fundamental recursion~(\ref{eq2}). Raising both members
of~(\ref{eq2}) to some integral power~$s$ yields 
\begin{equation}\label{hmom1}
X_n^s=X_{K_n}^s+\widetilde{X}_{n-a-K_n}^s+
\sum_{{s_1+s_2+s_3=s\atop s_2,s_3\not=s}}\binom{s}{s_1,s_2,s_3}
t_n^{s_1}X_{K_n}^{s_2}\widetilde{X}_{n-a-K_n}^{s_3},
\end{equation}
where we have  made use  of the multinomial expansion and
have isolated 
the two $s$th
powers. Take expectations with respect to the
model~$\mathfrak{M}_n$ and set~$\mu_n^{(s)}:=\E(X_n^s)$. The recursion
on $s$th moments becomes, thanks to independence  of the~$X$
and~$\widetilde{X}$ sequences on the right in~(\ref{hmom1}),
\begin{equation}\label{hmom2}
\mu_n^{(s)}=\sum_k p_{n,k}\left(\mu_k^{(s)}+\mu_{n-a-k}^{(s)}\right)+r_n^{(s)},
\end{equation}
where
\[
r_n^{(s)}:=\sum_{{s_1+s_2+s_3=s\atop
s_2,s_3\not=s}}\binom{s}{s_1,s_2,s_3}
t_n^{s_1}\sum_k p_{n,k}\mu_k^{(s_2)}\mu_{n-a-k}^{(s_3)}.
\]
This calculation shows that the sequence of~$s$th moments for any
fixed~$s$ satisfies the same type of recurrence as the first moments,
save  for a more complicated toll $(r_n^{(s)}$) that
involves moments of the smaller orders $0,1,\ldots,s-1$. Define the
normalized generating functions
\[
\mu^{(s)}(z):=\sum_n \mu_n^{(s)}\omega_n z^n,
\qquad
r^{(s)}(z):=\sum_n r_n^{(s)}\omega_n z^n,
\]
with the normalization sequence $\omega_n\equiv1$ for binary search
trees and $\omega_n=C_n$ for uniform binary  trees. Then
the relation~(\ref{hmom2}) is solved in terms of generating functions
by an~$\cal L$-transform as 
\begin{equation}\label{hmom3}
\mu^{(s)}(z)=\cal L\left(r^{(s)}\right),
\quad
r^{(s)}(z)=\!\!\!\!\!\!\!\sum_{{s_1+s_2+s_3=s\atop
s_2,s_3\not=s}}\binom{s}{s_1,s_2,s_3} \tau^{\odot s_1}(z) \odot
Q(\mu^{(s_2)}(z),\mu^{(s_3)}(z)),
\end{equation}
where
$Q$ is for the binary search tree model and Catalan model,
respectively, 
\begin{equation}\label{hmom4}
Q^{\operatorname{BST}}(a(z),b(z))=\int_0^z a(t) b(t)\, dt,
\qquad
Q^{\operatorname{Cat}}(a(z),b(z))=za(z)b(z).
\end{equation}
The~$\cal L$ transform is given in~(\ref{Lbst}) and~(\ref{Lcat}) for
the respective  cases; the case of the union--find tree model 
[where~$\omega_n=n^{n-1}/n!$ is used for $\mu^{(s)}$ and 
$\omega_n' = n^{n-2}(n-1)/n!$ is used for $r^{(s)}$] is similar
but more complicated---see~(\ref{Lcay}) for~$\cal L$, 
while
\[
Q^{UF}(a(z),b(z))=\frac12 a(z)b(z).
\]
As seen in the previous section, these $\cal L$ transforms involve
only  integration, differentiation, and  ordinary
and Hadamard products---all are operations that preserve the character
of being $\Delta$-regular and admitting complete asymptotic expansions
at the dominant singularity. We then have a general result:

\begin{theorem}\label{meta-thm} For any  of the binary search tree,
uniform binary tree, or union--find model, and for any integer~$s>0$, the
$s$th moment of the cost function associated to a toll $t_n^{\log}$ or
$t_n^\alpha$ admits a complete descending expansion in powers of~$n$ 
(possibly with logarithmic terms).
\end{theorem}
\begin{proof}
The proof is simply an induction on the order~$s$
of the moments. We establish by induction 
the stronger property that the generating functions 
$\mu^{(s)}(z)$ are $\Delta$-regular and admit complete asymptotic
expansions in powers of~$(1-z)$, possibly with logarithmic terms,
after rescaling the singularity to be at~1.
The property is true for~$s=1$
by results of the previous section. If the property is assumed to be
true through  order~$s-1$, then the tolls $r^{(s)}(z)$ 
are $\Delta$-regular and admit of complete
asymptotic expansions at their singularity: this results from
closure theorems
of Sections~\ref{sec:sing-expans-diff} and~\ref{sec:hadam-prod-transf}. Next,
the $\cal L$-transform is applied and, again by closure theorems, the property
of~$r^{(s)}(z)$
is seen to extend to $\mu^{(s)}(z)$. Thus the singular structure
of~$\mu^{(s)}(z)$  is fully characterized.
It then suffices to apply basic
singularity analysis in order to recover 
the existence of
full asymptotic expansion of the moments
$\mu_n^{(s)}=\frac{1}{\omega_n}[z^n]\mu^{(s)}(z)$. 
\end{proof}

The process of extracting moments one after the other has been
nicknamed ``moment pumping'' in the article~\cite{FlPoVi98}, where it
was used to determine the shape of the moments of total displacement
in linear hashing tables. It had been employed earlier by Louchard and
Tak\'acs in order to characterize moments of path length in trees and
of area under excursions~\cite{Louchard84,Takacs91}, in a way largely
 similar to what has been described here in more general
terms.  In favorable  cases, a pattern regarding the asymptotic shape
of moments may emerge. In such cases (possibly centering of the random
variable is required), the limiting distribution of costs becomes
accessible through its  moments, thanks to the moment convergence
theorem.  Instances are found in the already cited
papers~\cite{FlPoVi98,Louchard84,Takacs91}.  Fill's
paper~\cite{Fill96} provides another example (although it is based on
direct recurrence manipulations rather than generating functions) to
the effect that the logarithmic toll~$t^{\log}_n$ gives rise to
asymptotically Gaussian costs under the binary search tree model.  Yet
other examples, often based on direct recurrence manipulations, are
provided by the recent independent studies  of Hwang and
Neininger~\cite{HwNe02} and of Fill and 
Kapur~\cite{FK-transfer}. Clearly, a
``metatheorem'' similar
to Theorem~\ref{meta-thm} is possible for varieties of increasing
trees in the sense of Bergeron-Flajolet-Salvy~\cite{BeFlSa92}
(generalizing the BST model).  
For simply generated  families of trees in the
sense of Meir and Moon~\cite{MeMo78} (generalizing the Catalan model),
asymptotics of moments as well as limiting distributions have been
derived by Fill and Kapur~\cite{sgtechreport} as part of a
broader  project joint with Svante Janson.  The union--find
tree model can be generalized to other families of trees, and the
techniques of the present paper can again be applied; this is the
subject of ongoing research by the authors.

\subsection{Differential models}
Many tree recurrences
associated to comparison-based searching 
and multidimensional retrieval problems
generalizing binary search trees, once translated into generating
functions, lead to integral equations of the form
\begin{equation}\label{dif0}
{\bf \Phi}[f](z)=t(z),
\end{equation}
where ${\bf \Phi}$ is a linear integral operator involving 
coefficients in $\C(z)$, that is, rational function
coefficients. Here, as  in our lead examples, $f(z)$
is a generating function of expected costs and $t(z)$
is a toll generating function.
By successive differentiations, this transforms into a
linear differential equation of the form
\begin{equation}\label{dif1}
{\bf \Delta}[f](z)=\widetilde t(z), 
\end{equation}
where $\widetilde t(z)$ is a modified toll generating
function and  is an elementary variant of~$t(z)$.
We shall let~$d$ denote the order of the differential  equation~(\ref{dif1}).

The description above
corresponds to the situation already encountered with the binary search tree
recurrence, 
representing the easy case of a differential order equal to~1.
Other known cases include the $m$--ary search tree studied by Mahmoud
and Pittel (see the account in~\cite{Mahmoud92}) and 
others~\cite{ChHw01,FK-transfer}, 
quicksort with median-of-sample partitioning
and locally balanced trees~\cite{Knuth98a,SeFl96},
quadtrees~\cite{FlGoPuRo93,HoFl92}  as well as
multidimensional search trees also known as $k$-d-trees~\cite{FlPu86}.
(A valuable
survey of a class  of problems leading to Euler equations
appears in~\cite{ChHwTs02}.) 
For instance, in the case of 2-dimensional quadtrees  the operator
is given in~\cite{HoFl92} as 
\[
{\bf \Phi}[f](z)=f(z) - 4\int_0^z\left[\int_0^y f(x)\,
\frac{dx}{1-x}\right]\,\frac{dy}{y(1-y)},
\]
which leads to a second order differential equation, 
\[
z(1-z)\partial_z^2f(z)+(1-2z)\partial_zf(z)-\frac{4}{1-z} f(z) =
\widetilde t(z),
\]
where  $\widetilde t(z)=\partial_z[z(1-z)t'(z)]$.

The variation-of-constants technique applies to equations of order
greater than~1 as well as to linear systems.  It  may then be used to
express~$f$ as a linear integral transform involving a set $\{h_j\}$
of solutions to the homogeneous equation ${\bf \Delta} h=0$, as we
know explain.  Indeed, let the linear differential
equation~(\ref{dif1}) be put into the form of a system
\begin{equation}\label{difsys}
\partial_z {\bf y}(z) = {\bf A} {\bf y}(z)+ {\bf b}(z),
\end{equation}
where ${\bf y}$ is the $d$--dimensional vector ${\bf
  y}=(f,f',\ldots,f^{(d-1)})$, ${\bf A}={\bf A}(z)$ is a $d\times d$
matrix of functions  [here, by assumption, all in~$\C(z)$],
and~${\bf b}=(\widetilde{t},\widetilde{t}\,',\ldots,
\widetilde{t}^{(d-1)})$;  see~\cite[vol~II,
\S9.3]{Henrici77a} for the reduction.  Recall that a fundamental
matrix ${\bf W}$ for the system~(\ref{difsys}) is by definition a
nonsingular~$d\times d$ matrix whose columns each satisfy the
homogenous system $\partial_z {\bf y}(z) = {\bf A} {\bf y}(z)$.  Then
 the general solution to the inhomogeneous
system~(\ref{difsys}) is, by the classical ``variation-of-constants''
formula,
\begin{equation}\label{vofc}
{\bf y}(z)={\bf W}(z)\cdot{\bf W}^{-1}(z_0)\cdot {\bf y}(z_0)
+{\bf W}(z)\cdot \int_{z_0}^z {\bf W}^{-1}(x)\cdot  {\bf b}(x)\, dx;
\end{equation}
see once more~\cite[vol~II]{Henrici77a}. 
(The initial conditions at some~$z_0$ are assumed to be known.)
This provides the  solution to~(\ref{dif0}) as
\[
f(z)={\cal L}[\widetilde t(z)], 
\]
with~$\cal L$ a \emph{linear integral transform} that involves 
polynomially the elements of a fundamental matrix~${\bf W}$ as well as the
inverse of the Wronskian $\det {\bf W}$.
(The case of Euler equations
is somewhat simpler,  as it is fully
explicit~\cite{ChHwTs02,FK-transfer}.)  
For instance, the case of 2-dimensional quadtrees leads to
a still explicit form~\cite{HoFl92}, namely,
$
f(z)={\cal L}[\widetilde t(z)]$ where
\[
{\cal L}[e(z)]=\frac{1+2z}{(1-z)^2}
\int_0^z \frac{(1-y)^3}{y(1+2y)^2}\left[
\int_0^y \frac{1+2x}{(1-x)^2} e(x)\, dx\right]\, dy.
\]

From here onward we suppose for simplicity that $\omega_n \equiv 1$,
so that $f$ is an ordinary generating function, though our discussion
extends readily to more general normalization constants.  Call a
system \emph{dominantly regular}\footnote{The term ``dominantly
  regular'' evokes the fact that the condition concerns the
  \emph{dominant} singularity of the solution function, where the
  singularity of the system is of the so-called ``\emph{regular}''
  type (first-kind implies regular singularity by a well-known
  theorem; see~\cite[vol~II, Theorem~9.4d]{Henrici77a}).} 
if it is singular at~1 (i.e., if the matrix ${\bf A}$ has a pole at~1,
but at no other point in~$|z|\le1$ except possibly~0) and if the pole
of~${\bf A}$ at~1 is simple---the latter case is known as a
singularity ``of the first kind''.  All the classical examples listed
above and generalizing binary search trees satisfy this condition.  We
then have:
\begin{theorem}\label{genBST-thm}
Let a tree recurrence be expressed by a differential system that is
dominantly regular. Then the expectations of
costs induced by the tolls~$t_n^\alpha$
and~$t_n^{\log}$ admit complete asymptotic expansions in
descending powers of~$n$, possibly with logarithmic terms. 
\end{theorem}
\begin{proof}
  First, we observe that for any tree recurrence, the cost induced by
  an eventually increasing nonnegative toll $t_n \to +\infty$ is at least~$t_n$
  (by the very nature of the tree recurrence) and at most $O(nt_n)$
  (by induction). Thus, for the tolls under consideration, the
  generating function of costs, $f(z)$, has radius of convergence
  exactly equal to~1. We also observe that the values of~$f$ and its
  derivatives at some point~$z_0$ such that $|z_0|<1$ are well-defined.
  We may adopt for
  instance~$z_0=\frac12$ in the variation-of-constants
  formula. 
  
  By the classical theory of singularities of the first kind, each of
  the column vectors of matrix~${\bf W}$ is analytic for $z$ in a
  neighborhood
  of~1 slit along the ray $[1,+\infty)$.
There, as $z$ tends to~1, it admits a representation
as a finite combination of terms of the form
\[
(1-z)^{\alpha}L(z)^k R(1-z),
\]
where~$\alpha$ is an algebraic number (a root of the indicial
equation), $k$ an integer, and $R$ is analytic at~0. Thus, each
element of~${\bf W}$ is amenable to singularity analysis
as it is $\Delta$-continuable and admits a
   bona fide expansion near~1.
  
  By formula~(\ref{vofc}), there remains to discuss the elements of
  ${\bf W}^{-1}$. By the cofactor rule, the elements of ${\bf W}^{-1}$
  involve polynomially the elements of~${\bf W}$  divided by
  the Wronskian determinant $\det {\bf W}(z)$. It is a well-known fact
  (see {\S9.3} of~\cite[vol~II]{Henrici77a}) that the Wronskian is
  expressible in terms of the system alone and one has 
\[
[\det {\bf W}(z)]^{-1}= [\det {\bf W}(z_0)]^{-1}\exp
\left(-\int_{z_0}^z \operatorname{tr} {\bf A}(x)\, dx\right).
\]
[Here~$\operatorname{tr}(\cdot)$ denotes the matrix trace operator.] 
By the dominant regularity assumption, the trace is here a rational function
with at most a simple pole at~1, so that its integral is either
analytic at~1 or logarithmic. In either case, $(\det {\bf W}(z))^{-1}$
is of singularity analysis type, and so are consequently all the
elements of the inverse of the fundamental matrix~${\bf W}^{-1}(z)$.
By the singular integration and singular differentiation theorems of
Section~\ref{sec:sing-expans-diff}, there results that the integral
transform~(\ref{vofc}) preserves for functions the character of being
amenable to singularity analysis.  Since the toll generating functions
are of singularity analysis type, basic singularity analysis is
applicable to~$f(z)$.  The result follows.
\end{proof}

In principle, higher moments will also become accessible to
singularity analysis once the nonlinear integral forms~$Q$
extending $Q^{BST}$ of~(\ref{hmom4}) have been worked out.
We are however not aware of existing research in this direction, despite
the fact that the splitting probabilities are known in a number of cases
(see, e.g.,~\cite{FlGoPuRo93} for quadtrees). There
is interest in these questions,  as partly heuristic recent work by
Majumdar and collaborators (see, e.g.,~\cite{DeMa02} for the type of
method employed and succinct developments) indicates the 
probable existence of phase transitions in
the number of internal nodes of~$d$-dimensional quadtrees 
for large enough~$d$ ($d\ge d_c=9$ is suggested) in a way similar to what is 
already well established for the size of~$m$-ary search
trees~\cite{ChHw01,FK-transfer,Mahmoud92}.

\smallskip

As a final note, we'd like to mention digital trees, which were
recognized to be amenable to treatment by \emph{ordinary}
(rather than the more customary exponential) generating functions
in~\cite{FlRi92}. Techniques of the present paper 
would most likely be usable in such a context, in particular as
regards tolls of the form~$n^\alpha$ and~$\log n$. A partial
classification of cost functions along these lines has 
already been given by Derfel and Vogl in~\cite{DeVo01}.

\medskip
\noindent
\begin{small}
  {\bf Acknowledgements.} The second author is grateful to Don Knuth
  for a summer invitation to 
  Stanford in 1988, during which his work on the subject was started.
  This work was supported in part by the Future and Emerging
  Technologies programme of the EU under contract number
  IST-1999-14186 ({\sc Alcom--Ft} Project).

The first and third authors' research is supported by NSF grant
DMS--9803780 and DMS--0104167, and by The Johns Hopkins University's
Acheson J.~Duncan Fund for the Advancement of Research in Statistics.
\end{small}

\small

\bibliographystyle{acm}
\bibliography{algo,msn,misc}

\end{document}